\documentclass[12pt]{amsart}
\usepackage{amsfonts,amsmath,amsthm,amssymb}
\usepackage{psfrag}
\usepackage{graphicx}
\usepackage{epstopdf}
\usepackage{stmaryrd}
\usepackage[mathscr]{eucal}
\begin{document}
\newtheorem{theorem}{Theorem}
\newtheorem{proposition}{Proposition}
\newtheorem{lemma}{Lemma}
\newtheorem{definition}{Definition}
\newtheorem{notation}{Notation}
\newtheorem{remark}{Remark}

\title{\bf A Presentation For The Automorphisms Of The $3$-Sphere That  Preserve A Genus Two Heegaard Splitting}\author{Erol Akbas}\maketitle
\date{}

\begin{abstract}
  {\it  Scharlemann constructed a connected simplicial 2-complex  $\Gamma$ with an action by the group ${\mathcal H_{2}}$ of isotopy classes of orientation preserving homeomorphisms of $S^3$ that preserve the isotopy class of an unknotted genus $2$ handlebody $V$. In this paper we prove that the 2-complex $\Gamma$ is contractible. Therefore we get a finite presentation of  ${\mathcal H_{2}}$.}
\end{abstract}

\section{Introduction}
Let ${\mathcal H_{g}}$ be the group of isotopy classes of orientation preserving homeomorphisms of $S^3$ that preserve the isotopy class of an unknotted genus $g$ handlebody $V$. In $1933$ Goeritz \cite{Go} proved that $\mathcal{H}_2$ is finitely generated. In $1977$ Goeritz's theorem was generalized to arbitrary genus $g\geq 2$ by Jerome Powell \cite{Po}.  In $2003$ Martin Scharlemann   noticed that Powell's proof contains a serious gap. Scharlemann \cite{Sc} gave a modern proof of Goeritz's theorem by introducing a  simplicial $2$-complex $\Gamma$, with an action by ${\mathcal H_{2}}$, that deformation retracts onto a graph $\tilde{\Gamma}$. Given any two distinct vertices  $v$, $\tilde{v}$ of $\Gamma$, Scharlemann constructed a vertex $u$  in $\Gamma$ that is adjacent to $v$ and ``closer" to $\tilde{v}$  (by ``closer" we mean the intersection number of $u$ and $\tilde{v}$, see Definition \ref{defn:intersectionnumber}). Hence ${\mathcal H_{2}}$ acts on the connected graph $\tilde{\Gamma}$ and is generated by the isotopy classes of elements denoted by $\alpha$, $\beta$, $\gamma$ and $\delta$ (see Section \ref{sec:preliminaries} for a complete description). In this paper we  study the geometry of $\Gamma$ by showing that $u$ is essentially unique (for a precise statement see Proposition  \ref{backbone}). We derive the following theorem.

\begin{theorem}\label{tree}
The graph $\tilde{\Gamma}$ is a tree, and shortest paths can be calculated algorithmically.
\end{theorem}

Note that $\tilde{\Gamma}$ is locally infinite. So calculating paths is not trivial. We also get 

\begin{theorem}\label{presentation}
\begin{description}
\item[i] $\mathcal{H}_2=<[\alpha],[\beta],[\gamma],[\delta]\ |\ [\alpha]^2=[\gamma]^2=[\delta]^3=[\alpha\gamma]^2=[\alpha\delta\alpha\delta^{-1}]=[\alpha\beta\alpha\beta^{-1}]=1,[\gamma\beta\gamma]=[\alpha\beta],[\delta] =[\gamma\delta^2\gamma]>$

\item[ii]  $\mathcal{H}_2 \cong  (\mathbb{Z} \oplus \mathbb{Z}_2) \rtimes \mathbb{Z}_2 \underset{\mathbb{Z}_2 \oplus \mathbb{Z}_2}{\ast} (\mathbb{Z}_3 \rtimes \mathbb{Z}_2) \oplus \mathbb{Z}_2$
\end{description}
\end{theorem}

{\bf Acknowledgement}

This paper is a result of my thesis. I am grateful to my advisor, Professor Marc Culler, for his kind supervision and support. I would like to thank Ian Agol for his suggestions for the proofs of Lemma \ref{notequal} and Lemma \ref{lemmaC}. I would also like to thank Yo'av Rieck for his helpful conversations.

\section{Preliminaries}\label{sec:preliminaries}
We give a description of the  $2$-complex $\Gamma$ introduced by Scharlemann in \cite{Sc}. For details about $\Gamma$ we refer the reader to \cite{Sc}.

Let $V$ be an unknotted  handlebody of genus two in $S^3$ and let $W$ be the closure of its complement. Let $T$ be the boundary of $V$. Then $T$ is a genus two Heegaard surface for $S^3$. Let ${\mathcal H_{2}}$ denote the group of isotopy classes of orientation preserving homeomorphisms of $S^3$ that leave  the genus two handlebody $V$ invariant. A sphere $P$ in $S^3$ is called a {\em reducing sphere} for $T$ if $P$ intersects $T$ transversely in a simple closed curve which is homotopically non-trivial on $T$. For any reducing sphere $P$ for $T$ let $c_P$ denote $P\cap T$ and let $v_P$ denote the isotopy class of $c_P$ on $T$.

\begin{definition} \label{defn:intersectionnumber}
For any two reducing spheres $R$, $Q$ for $T$,  define the intersection number of $v_R$ and $v_Q$ as  
\begin{displaymath}
v_{R}\cdot v_{Q}=\underset{R^\prime\in v_R\ \ Q^\prime\in v_Q}{\min} |c_{R^\prime}
\cap c_{Q^\prime}|
\end{displaymath} 
where $|c_{R^\prime}\cap c_{Q^\prime}|$ is the geometric intersection number of $c_{R^\prime}$ with $c_{Q^\prime}$.
\end{definition}

Let $\Gamma$ be a complex whose vertices are isotopy classes of reducing spheres for $T$. A collection $P_0,...,P_n$ of reducing spheres bounds an $n$-simplex in $\Gamma$ if and only if $v_{P_i} \cdot v_{P_j}=4$ for all $0 \leq i \neq j \leq n$. In fact $n\leq 2$ \cite[Lemma 2.5]{St}. So $\Gamma$ is a simplicial $2$-complex. Let $\triangle$ be any $2$-simplex of $\Gamma$. We denote by $S_\triangle$ the ``spine" of  $\triangle$, which is the subcomplex of the barycentric subdivision consisting of all closed $1$-simplices that contain the barycenter and a vertex of $\triangle$. Clearly $\triangle$ deformation retracts onto $S_\triangle$. Let $\tilde{\Gamma}$ be
\begin{displaymath}
\underset{\triangle}{\bigcup}S_\triangle.
\end{displaymath} 
So $\tilde{\Gamma}$ is a graph. Since no two $2$-simplicies of $\Gamma$ share an edge \cite[Lemma 2.5]{St}, the simplicial 2-complex $\Gamma$ deformation retracts onto the graph $\tilde{\Gamma}$.

\begin{figure}[h]
\begin{center}
\psfragscanon
\psfrag{a}[][][0.8]{$A$}
\psfrag{b}[][][0.8]{$B$}
\psfrag{c}[][][0.8]{$C$}
\psfrag{x}[][][0.8]{$X$}
\psfrag{y}[][][0.8]{$Y$}
\psfrag{z}[][][0.8]{$Z$}
\psfrag{cp}[][][0.8]{$c_P$}
\psfrag{t-}[][][1.2]{$T^-$}
\psfrag{t+}[][][1.3]{$T^+$}
\includegraphics[height=3.5cm, width=5cm]{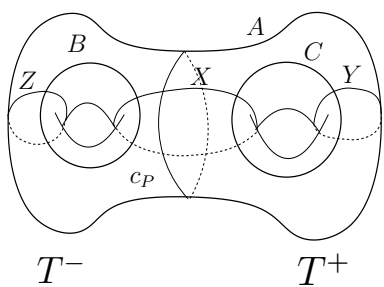}
\end{center}
\caption{The curves $c_P$, $A$, $B$, $C$, $X$, $Y$ and $Z$}\label{fig:abcxyz}
\end{figure}

A \textit{belt curve} on a genus two surface is a homotopically nontrivial separating simple closed curve. 
Let $P$ denote a reducing sphere whose intersection with $T$ is a belt curve, which we denote $c_{P}$.
The reducing sphere $P$ divides $S^3$ into two 3-balls $B^\pm$ whose intersections with the genus two surface $T$ are two genus one surfaces $T^\pm=T\cap B^\pm$, each having one boundary component. The surface $T^-$  (resp. $T^+$) contains two simple closed curves $B$, $Z$  (resp.  $C$, $Y$) meeting at one point. The curve  $B$ (resp. $C$) bounds a non-separating disc in $W$, homotopically non-trivial in $V$. The curve $Z$  (resp. $Y$) bounds a non-separating disc in $V$, homotopically non-trivial in $W$. The genus two surface $T$ contains two disjoint simple closed curves $A$ and $X$. The curve $A$ is homotopically non-trivial in $V$, disjoint from $B$ and $C$, bounds a non-separating disc in $W$, and intersects $Z$ and $Y$  at one point. The curve  $X$  is homotopically non-trivial in $W$, disjoint from $Z$, $Y$ and $A$, bounds a non-separating disc in $V$, and intersects $B$ and $C$ at one point. See figure \ref{fig:abcxyz}.

Throughout this paper, unless otherwise stated, whenever we choose a  reducing sphere $R$ for $T$ such that $v_R\neq v_P$ we will assume that the curve $c_R$ intersects  $c_P$, $B$, $C$, $Y$, $Z$ transversely and minimally, and intersects  $A$ transversely.

\begin{figure}[h]
\begin{center}
\psfragscanon
\psfrag{alpha}[][][1.5]{$\alpha$ :}
\psfrag{beta}[][][1.5]{$\beta$ :}
\psfrag{gamma}[][][1.5]{$\gamma$ :}
\psfrag{delta}[][][1.5]{$\delta$ :}
\psfrag{cp}[][][1]{$c_P$}
\psfrag{pi}[][][1.2]{$\pi$}
\psfrag{2pi/3}[][][1.2]{$\frac{2\pi}{3}$}
\psfrag{a}[][][1]{$A$}
\psfrag{b}[][][1]{$B$}
\psfrag{c}[][][1]{$C$}
\psfrag{x}[][][1]{$X$}
\psfrag{y}[][][1]{$Y$}
\psfrag{z}[][][1]{$Z$}
\includegraphics[height=2.5cm, width=12cm]{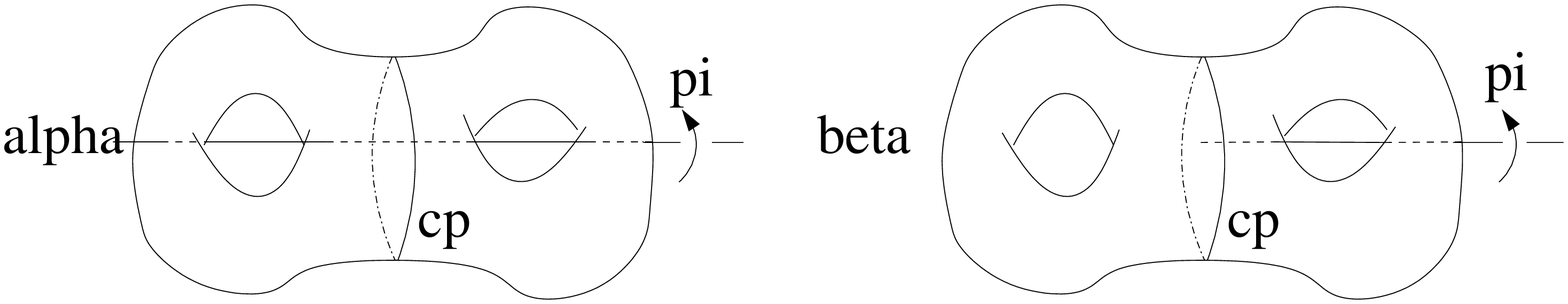}\vskip0.5cm
\includegraphics[height=3cm, width=6cm]{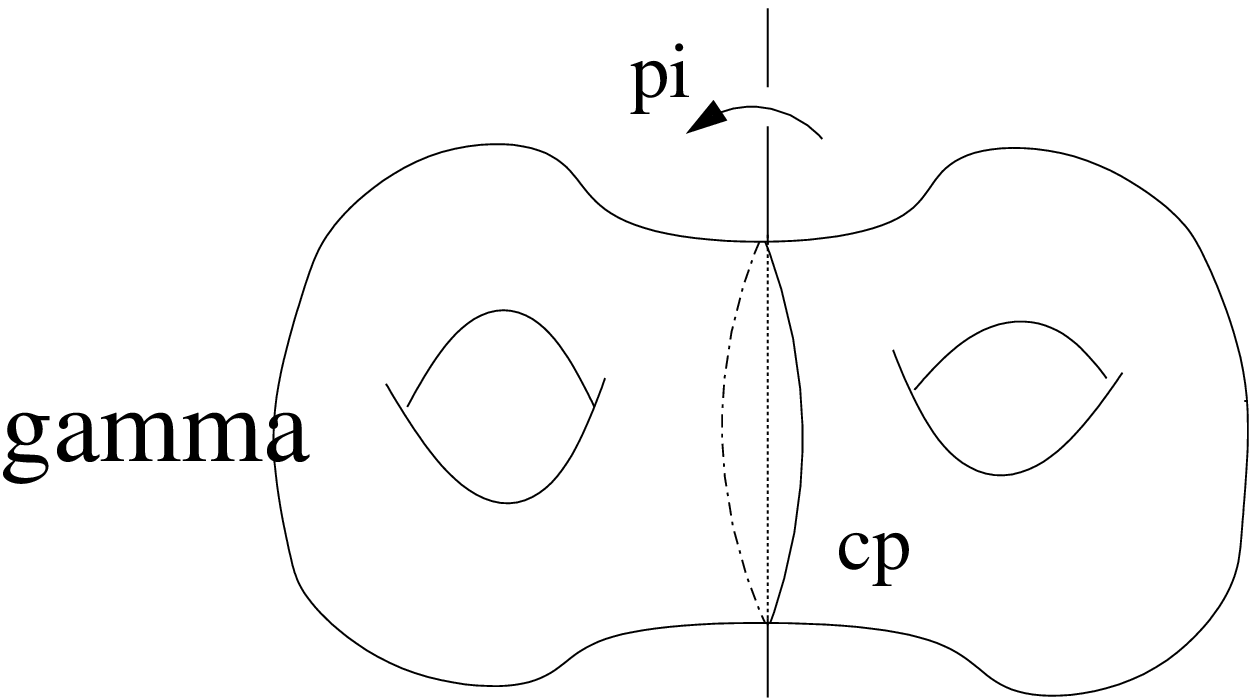}\hskip3cm
\includegraphics[height=4cm, width=3.5cm]{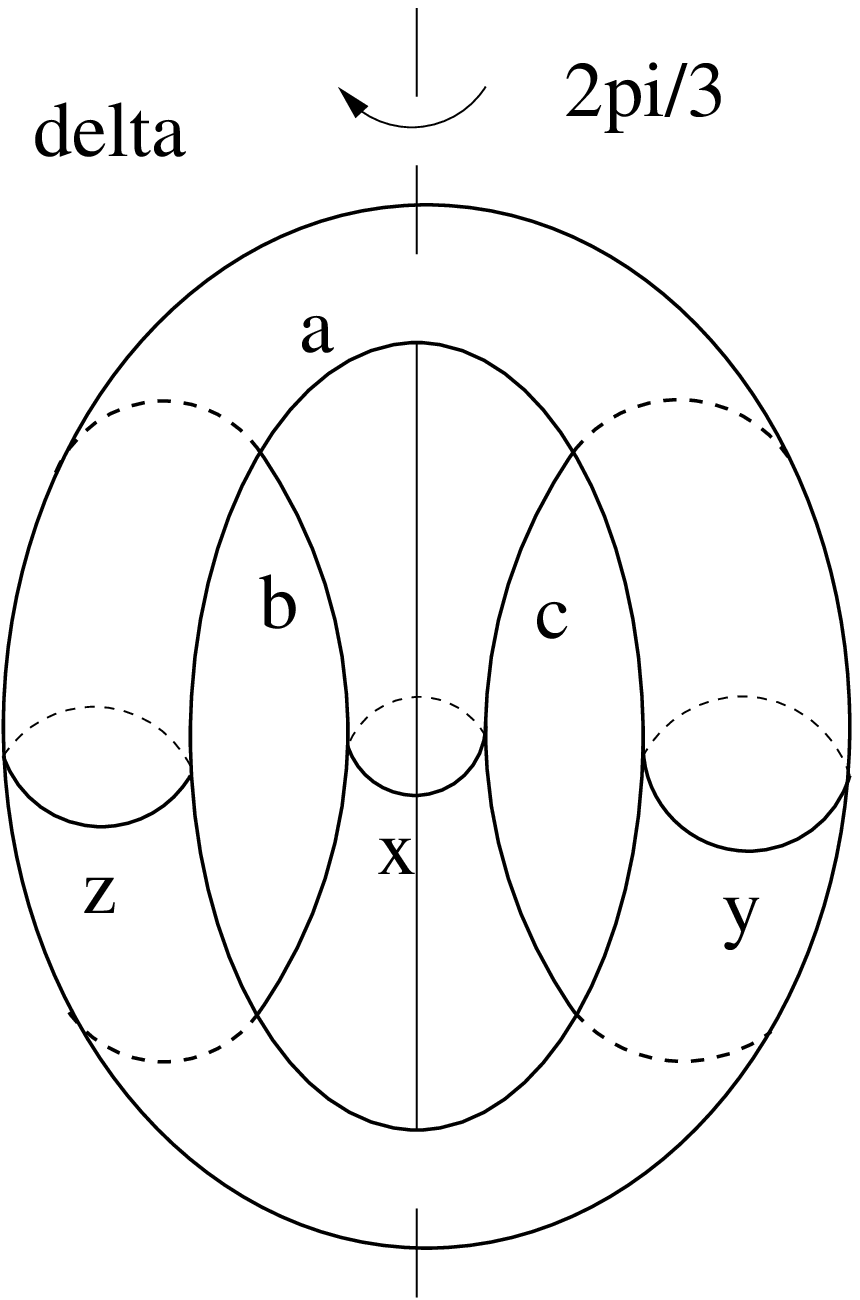}\hfill
\end{center}
\caption{Homeomorphisms $\alpha$, $\beta$, $\gamma$ and $\delta$}\label{fig:alphabetagammadelta}
\end{figure}

There exist three automorphisms $\alpha$, $\beta$, $\gamma$ of $S^3$  with the following properties. The automorphism $\alpha$  is an orientation preserving homeomorphism of $S^3$ that preserves  $V$ and $P$, and that maps the curves $A$, $B$, $C$ to  $A$, $B$, $C$ respectively by an orientation reversing map. The homeomorphism $\alpha$ is the hyperelliptic involution which preserves every simple closed curve (upto isotopy). The automorphism $\beta$ is an orientation preserving homeomorphism of $S^3$ that preserves  $V$ and $P$, fixes $T^-$ pointwise, and maps $C$ to $C$ and  $Y$ to $Y$ by an orientation reversing map. Also $|A\cap \beta(X)|=2$.  The automorphism $\gamma$  preserves  $V$ and $P$, and maps the curves $c_P$ to $c_P$ and $A$ to $A$ by an orientation reversing map. See figure \ref{fig:alphabetagammadelta}. Scharlemann \cite{Sc} showed that $\mathcal{H}_2$ is generated by the isotopy classes $[\alpha]$,  $[\beta]$,  $[\gamma]$ and $[\delta]$ where  $\delta$ is any orientation preserving homeomorphism of $S^3$ such that $\delta(V)=V$ and $v_{P}\cdot v_{\delta(P)}= 4$. In this paper we will take $\delta$ as the following homeomorphism. Consider the genus two handlebody $V$ as a regular neighborhood of a sphere, centered at the origin, with three holes. The homeomorphism $\delta$  is $\frac{2\pi}{3}$ rotation  of $V$ about the vertical $z$-axis. See figure \ref{fig:alphabetagammadelta}. 

\section{Arc families of reducing spheres on $T^\pm$}
\begin{definition}\label{def:slope}
 For any properly embedded arc $\nu\subset T^\pm$  we may write $[\nu]\in H_1(T^\pm,\partial T^\pm;\mathbb{Z})$ as $a[\mu]+b[\lambda]$ where $\mu=Z$ and  $\lambda=B$ if $\nu\subset T^-$, and $\mu=Y$ and $\lambda=C$ if  $\nu\subset T^+$. The slope of $\nu$ is defined to be  $|\frac{a}{b}|\in \mathbb{Q^+}\cup \infty$.
\end{definition}

\begin{definition}
For any reducing sphere $Q$  such that $v_Q\neq v_P$, let $N(Q,T^\pm,a)$ denote the number of arcs in $Q\cap T^\pm$ of slope $a$.
\end{definition}

\begin{definition}   
Denote any oriented curve $D$  on $T$ by $\overset{\rightarrow}{D}$ and the curve oriented in the  direction opposite to $\overset{\rightarrow}{D}$ by $\overset{\leftarrow}{D}$.
\end{definition}

\begin{definition}\label{defn:Theta}
Orient the curves $A$, $B$, $C$, $X$, $Y$, $Z$ in such a way that
$\delta^2(\overset{\rightarrow}{A})
=\delta(\overset{\rightarrow}{B})
=\overset{\rightarrow}{C}$ and
$\delta^2(\overset{\rightarrow}{X})
=\delta(\overset{\rightarrow}{Y})
=\overset{\rightarrow}{Z}$.
Up to isotopy there are natural homeomorphisms $\Omega,  \ \Psi : S^3 \longrightarrow  S^3$ where $\Omega$ maps $V$ to $W$ and
$\overset{\rightarrow}{A}$,
$\overset{\rightarrow}{B}$,
$\overset{\rightarrow}{C}$,
$\overset{\rightarrow}{X}$,
$\overset{\rightarrow}{Y}$,
$\overset{\rightarrow}{Z}$
to
$\overset{\leftarrow}{X}$,
$\overset{\leftarrow}{Y}$,
$\overset{\leftarrow}{Z}$,
$\overset{\rightarrow}{A}$,
$\overset{\rightarrow}{B}$,
$\overset{\rightarrow}{C}$ respectively,
and $\Psi$ maps $W$ to $W$ and
$\overset{\rightarrow}{A}$,
$\overset{\rightarrow}{B}$,
$\overset{\rightarrow}{C}$,
$\overset{\rightarrow}{X}$,
$\overset{\rightarrow}{Y}$,
$\overset{\rightarrow}{Z}$
to
$\overset{\rightarrow}{A}$,
$\overset{\rightarrow}{B}$,
$\overset{\rightarrow}{C}$,
$\overset{\leftarrow}{X}$,
$\overset{\leftarrow}{Y}$,
$\overset{\leftarrow}{Z}$ respectively (see figure \ref{fig:arcsonT}). Let $\Theta=\Psi\Omega$. 
\begin{figure}[h]
\begin{center}
\psfragscanon
\psfrag{1}[][][1]{$T^-$}
\psfrag{2}[][][1]{$T^+$}
\psfrag{4}[][][1.2]{$\Omega$}
\psfrag{5}[][][1.2]{$\Psi$}
\psfrag{ar}[][][0.8]{$\overrightarrow{A}$}
\psfrag{br}[][][0.8]{$\overrightarrow{B}$}
\psfrag{cr}[][][0.8]{$\overrightarrow{C}$}
\psfrag{xr}[][][0.8]{$\overrightarrow{X}$}
\psfrag{yr}[][][0.8]{$\overrightarrow{Y}$}
\psfrag{zr}[][][0.8]{$\overrightarrow{Z}$}
\psfrag{xl}[][][0.8]{$\overleftarrow{X}$}
\psfrag{yl}[][][0.8]{$\overleftarrow{Y}$}
\psfrag{zl}[][][0.8]{$\overleftarrow{Z}$}
\psfrag{cf}[][][0.8]{$c_f$}
\includegraphics[height=7cm, width=12cm]{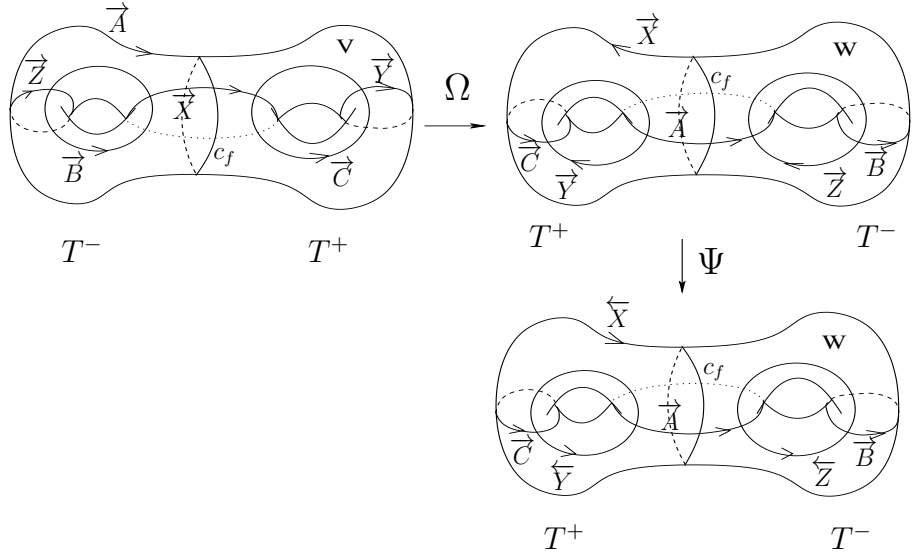}
\end{center}
\caption{Homeomorphism $\Theta=\Psi\Omega$}\label{fig:arcsonT}
\end{figure}
\end{definition}

\begin{proposition} \label{a=1/a}
Let $Q$ be a  reducing sphere for $T$ such that  $v_Q\neq v_P$. Then $N(Q,T^-,a)= N(Q,T^+,\frac{1}{a})$.
\end{proposition}

\begin{proof}
Without loss of generality, we may assume that $Q=w(P)$ where $w$ is a word in $\alpha$, $\beta$, $\gamma$ and $\delta$. 

\textbf{Claim.} $\Theta(c_Q)=c_Q$

\textbf{Proof of Claim.}  The hyperelliptic involution $\alpha$ preserves the isotopy class of any simple closed curve on $T$. After an isotopy, we may assume that $\alpha(c_Q)=c_Q$.  Let us write $w$ as $a_1a_2...a_n$ where $a_i\in\{\alpha,\beta^{\pm1},\gamma,\delta^{\pm1}\}$. The homeomorphism  $\Theta$ satisfies $\Theta\alpha=\alpha\Theta$, $\Theta\beta=\alpha\beta\Theta$, $\Theta\gamma=\alpha\gamma\Theta$, $\Theta\delta=\delta\Theta$, and  $\Theta(c_P)=c_P$. Then $\Theta(c_Q)=\Theta(w(c_P))=\Theta(a_1a_2...a_n(c_P))=b_1b_2...b_n\Theta(c_P)$ where $b_i$ is $\alpha$ if $a_i=\alpha$, $\alpha\beta$ if $a_i=\beta$, $\alpha\gamma$ if $a_i=\gamma$, $\delta$ if $a_i=\delta$. So $b_1b_2...b_n\Theta(c_P)=b_1b_2...b_n(c_P) = a_1a_2...a_n(c_P)=w(c_P)=c_Q$.

Since $\Theta$ maps the curves  $A$, $B$, $C$, $X$, $Y$, $Z$ to $X$, $Y$, $Z$, $A$, $B$, $C$ respectively, it  takes the arcs of $c_Q$ of slope  $a$ on $T^-$ to the arcs of $c_Q$ of slope  $\frac{1}{a}$ on $T^+$.
\end{proof}

\begin{lemma} \label{notequal}
Let $Q$ be any reducing sphere for $T$ such that $v_Q\neq v_P $.  Then $N(Q,T^-,0)\neq N(Q,T^-,\infty)$.
\end{lemma}

\begin{proof}
Suppose that $N(Q,T^-,0)=N(Q,T^-,\infty)=m$. By Proposition \ref{a=1/a}, $N(Q,T^+,0)=N(Q,T^+,\infty)=m$ and $N(Q,T^-,1)=N(Q,T^+,1)$. The curve $c_Q$ bounds a disc in $V$. So $c_Q$ must have a ``wave" $\tau$ \cite{vkf} with respect to one of the curves $Y$, $Z$. Say with respect to $Y$. Then the arc $\tau$ of $c_Q$  starts at $Y$, goes to $T^-$ then comes back to $Y$ on the same side without touching $Z$. So all the arcs of $c_Q$ intersecting $Z$ must  intersect the arc on $Y$ that is bounded by ends of $\tau$.  Then we get $N(Q,T^-,\infty)+N(Q,T^-,1)+2\leq N(Q,T^+,\infty)+N(Q,T^+,1)$, a contradiction.
\end{proof}

\begin{definition}
For any reducing sphere $Q$ for $T$ such that $v_Q\neq v_P$, let $F^\pm_{Q,a}$ denote  the arc family of $c_Q$ on $T^\pm$ of slope $a$.
\end{definition}

\begin{notation}\label{ends}
We will fix the following notation: Let $Q$ be a reducing sphere for $T$.

\begin{itemize}
\item If  $N(Q,T^-,0)= n\neq 0$ then $e_{01}$, $e_{02}$,..., $e_{0n}$, $e_{1n}$, $e_{1n-1}$,..., $e_{11}$ are going to denote consecutive end points, on $c_P$, of the arcs in  $F^-_{Q,0}$   where $e_{0j}$, $e_{1j}$  are end points of the same arc, and  $h_{01}$, $h_{02}$,..., $h_{0n}$, $h_{1n}$, $h_{1n-1}$,..., $h_{11}$ are going to denote consecutive end points, on $c_P$, of the arcs in $F^+_{Q,\infty}$  where $h_{0j}$, $h_{1j}$ are end points of the same arc (existence of $h_{ij}$ is guaranteed by Proposition \ref{a=1/a}).

\item If  $N(Q,T^-,\infty)= m\neq 0$ then  $g_{01}$, $g_{02}$,..., $g_{0m}$, $g_{1m}$, $g_{1m-1}$,..., $g_{11}$ are going to denote consecutive  end points, on $c_P$, of the arcs in  $F^-_{Q,\infty}$  where $g_{0j}$, $g_{1j}$  are end points of the same arc,  and  $f_{01}$, $f_{02}$,..., $f_{0m}$, $f_{1m}$, $f_{1m-1}$,..., $f_{11}$ are going to denote consecutive  end points, on $c_P$,   of the arcs in  $F^+_{Q,0}$ where $f_{0j}$, $f_{1j}$  are end points of the same arc.

\item If  $N(Q,T^-,1)= p \neq 0$  then  $k_{01}$, $k_{02}$,..., $k_{0p}$, $k_{1p}$, $k_{1p-1}$,..., $k_{11}$ are going to denote consecutive  end  points, on $c_P$, of the arcs in $F^-_{Q,1}$ where $k_{0j}$, $k_{1j}$  are end points of the same arc,  and  $l_{01}$, $l_{02}$,..., $l_{0p}$, $l_{1p}$, $l_{1p-1}$,..., $l_{11}$ are going to denote  end points, on $c_P$,  of the  arcs in  $F^+_{Q,1}$ where $l_{0j}$, $l_{1j}$ are  end points  of the same arc.
\end{itemize}
\end{notation}

\begin{lemma} \label{lemmaC}
Let $Q$ be a reducing sphere for $T$ such that  $N(Q,T^-,0) = n > N(Q,T^-,\infty) = m > N(Q,T^-,1) = 0$.  Then
$\{f_{ij}|i=0,1 \ \ j=1,m\}\subseteq\{e_{ij}|i=0,1 \ \ j=2,...,n-1\}$.
\end{lemma}

\begin{proof} Suppose that $\{f_{ij}|i=0,1 \ \ j=1,m\}\nsubseteq\{e_{ij}|i=0,1 \ \ j=2,...,n-1\}$ (see figure \ref{fig:lemmaC}). 
\begin{figure}[h]
\begin{center}
\psfragscanon
\psfrag{t-}[][][1.3]{$T^-$}
\psfrag{t+}[][][1.3]{$T^+$}
\psfrag{n}[][][1]{$n$}\psfrag{m}[][][1]{$m$}
\includegraphics[height=5.5cm, width=5cm]{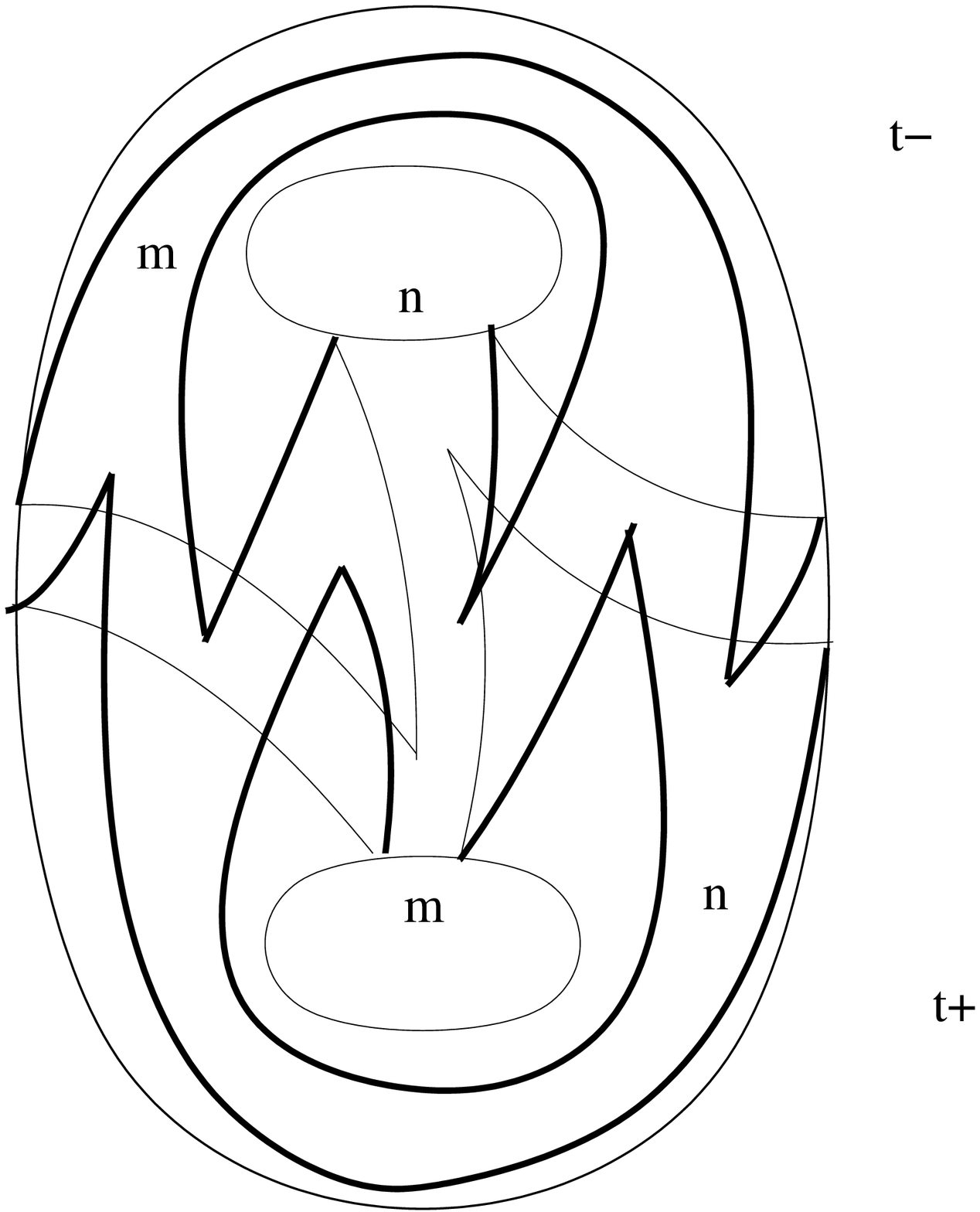}
\end{center}
\caption{}\label{fig:lemmaC}
\end{figure}
Then $c_Q$ does not have a ``wave" $\tau$ \cite{vkf} with respect to the curve $Y$ or the curve $Z$. Therefore $c_Q$ can not bound a disc in $V$, a contradiction.
\end{proof}

\newpage
\begin{proposition}\label{backbone}
Let $v$, $\tilde{v}$ be any two distinct vertices of $\Gamma$ such that $v\cdot\tilde{v}\neq4$. Then there exists unique vertex $u$ of $\Gamma$ such that
\begin{description}
\item[i] $u  \cdot v=4$
\item[ii] $u  \cdot \tilde{v} < v\cdot \tilde{v}$
\item[iii] $u  \cdot \tilde{v}  < v^\prime \cdot \tilde{v}$ for any vertex $v^\prime$ of $\Gamma$ such that $v^\prime \neq u$ and $v^\prime \cdot v=4$
\end{description}
Moreover, there is at most one vertex $v^{\prime\prime}$ of $\Gamma$ satisfying  $v \cdot v^{\prime\prime} =4$ and $u \cdot \tilde{v} < v \cdot v^{\prime\prime} \leq v \cdot \tilde{v}$. In this case $v^{\prime\prime} \cdot u=4$.
\end{proposition}

\begin{proof}
Let $v$, $\tilde{v}$ be any two vertices of $\Gamma$ such that $v\neq \tilde{v}$ and $v \cdot \tilde{v}\neq 4$. Since the group $\mathcal{H}_2$ is transitive on the vertices of $\Gamma$, we may assume that $v=v_P$ and $\tilde{v}$ is a vertex of $\Gamma$ such that $\tilde{v}\neq v_P$ and $v_P \cdot \tilde{v}\neq 4$. Then  for some word $w$ in $\alpha,\gamma,\beta$ and $\delta$, \ $w(P)\in \tilde{v}$.  Let $Q$ denote the reducing sphere $w(P)$.  Since $Q$ is not isotopic to $P$ there must be some arcs  in $c_Q\cap T^\pm$. By  \cite[Lemma 4]{Sc}  there is an arc of $c_Q$ of slope $0$  either on  $T^-$ or  on $T^+$.  Suppose it is on $T^-$.   Let $e_{ij}$, $g_{dq}$, $k_{rs}$, $f_{tu}$, $h_{yv}$, $l_{wz}$ denote the end points of the arcs of $c_Q\cap T^\pm$ as in the  Notation \ref{ends}. Possible cases for the arc families in  $c_Q\cap T^\pm$ and their configurations, upto a power of $\beta$, are the following:\\

\noindent {\bf Case I.} If $N(Q,T^-,0)=m$,  $N(Q,T^-,1/k)=a$ and  $N(Q,T^-,1/(k+1))=b$ where $k\geq 1$  then $N(Q,T^+,\infty)=m$,  $N(Q,T^+,k)=a$ and  $N(Q,T^+,k+1)=b$ by Proposition \ref{a=1/a}. Scharlemann in \cite[Lemma 5]{Sc}  constructs a reducing sphere $R$ satisfying (i) and (ii) (i.e. $v_R \cdot v_P = 4$  and $v_R \cdot  v_Q < v_P\cdot v_Q$). We will show that  upto isotopy the reducing sphere $R$ also satisfies (iii). Let $n=a+b$.

\noindent {\bf I.A.} If $n\neq 0$: Let us label end points of the arcs in $c_Q\cap T^+$ of slope  different from  $\infty$ as $d_{1}$, $d_{2}$,..., $d_{2n}$. Then it is not hard to show $\{e_{ij}\}\nsubseteq \{d_{i}\}$ by an argument similar to the proof of Lemma \ref{lemmaC}.

\noindent {\bf I.A.1.}  If $\{e_{ij}\}\cap\{h_{ij}\}\neq \emptyset$ (see figure \ref{fig:IA2}): Set $p=\frac{|\{e_{ij}\}\cap \{h_{ij}\}|}{2}$ then $1\leq p<m$. Consider the curve $\xi$ shown in figure \ref{fig:IA2}. It is easy to see that $\xi$ bounds a disc in $V$ and a disc in $W$. So $\xi$ is the intersection of a reducing sphere $S$ with $T$. Denote $\xi$ by $c_S$. The reducing sphere $S$ satisfies $v_S\cdot v_Q\leq |c_S\cap c_Q|=2(n-m+2p)<2(n+m)=v_P\cdot v_Q$ and $v_S\cdot v_P=4$.
\begin{figure}[h]
\begin{center}
\psfrag{f0-}[][][1]{$F^-_{Q,0}$}
\psfrag{f1k-}[][][1]{$F^-_{Q,\frac{1}{k}}$}
\psfrag{f1k+1-}[][][1]{$F^-_{Q,\frac{1}{k+1}}$}
\psfrag{fi+}[][][1]{$F^+_{Q,\infty}$}
\psfrag{cr}[][][1.5]{$\xi$}
\psfrag{mix1}{$F^-_{Q,\frac{1}{k}}\cup F^-_{Q,\frac{1}{k+1}}$}
\psfrag{mix2}{$F^+_{Q,k}\cup F^+_{Q,k+1}$}
\includegraphics[height=6.5cm, width=9cm]{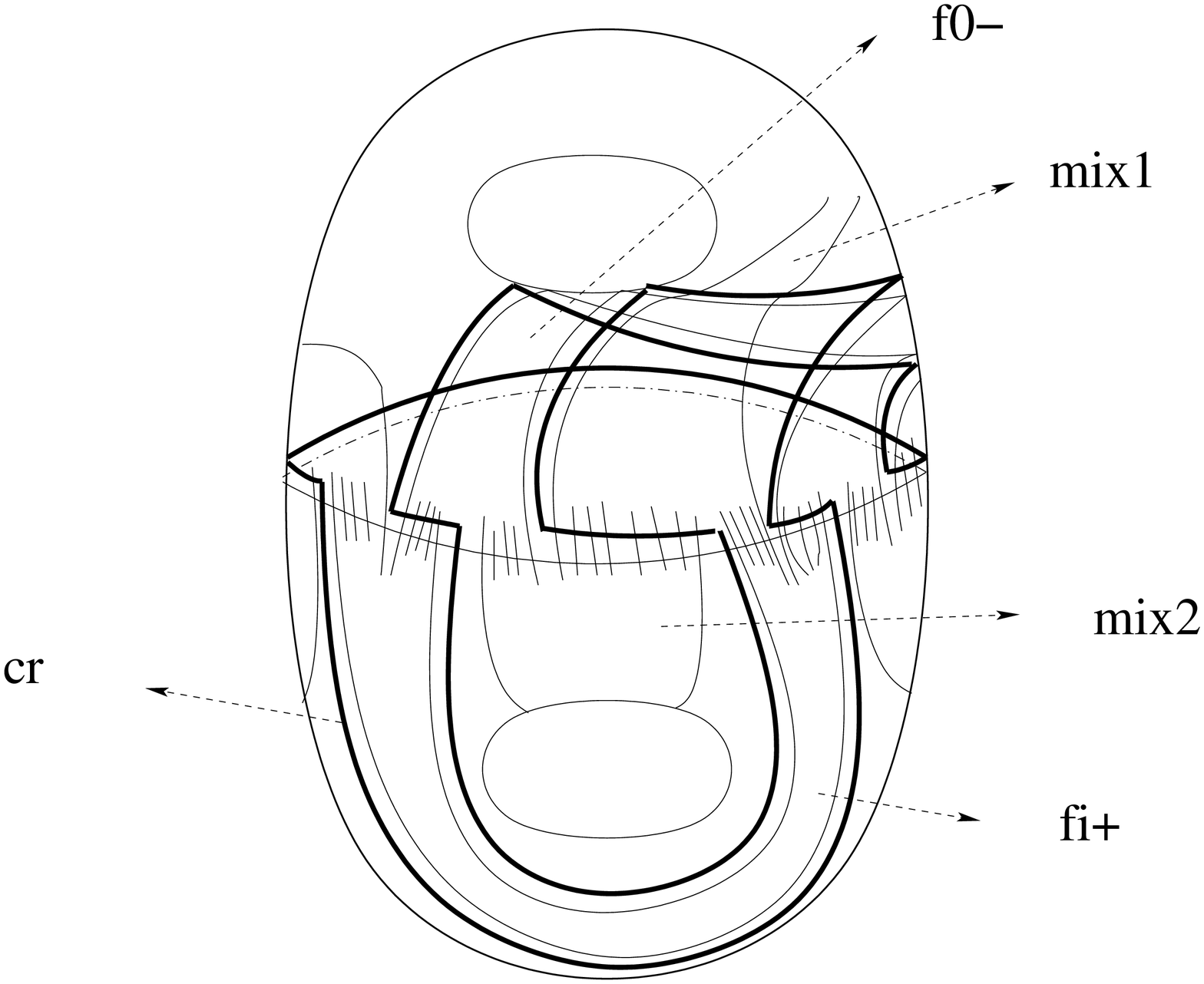}
\end{center}
\caption{}\label{fig:IA2}
\end{figure}

\textbf{Claim 1.} $v_S\cdot v_Q=|c_S\cap c_Q|$.

\textbf{Claim 2.} $v_{\beta^i (S)} \cdot v_Q,\ \ v_{\beta^i\gamma(S)}\cdot v_Q > 2(n+m)$ for $i\neq 0$.

\textbf{Proof of Claim 1.} It suffices to show that there is no bigon on $T$ formed by the curves $c_S$ and $c_Q$. We may assume that $c_S$ intersects $c_Q$ in a neighborhood $N\subseteq T$ of $c_P$ where $N\cap (B\cup Z\cup C\cup Y)=\emptyset$. The neighborhood $N$ has two boundary components $N^-$, $N^+$. Say $N^\pm\subset T^\pm$. The set $c_S\cap N$ consists of four arcs $\nu_1$, $\nu_2$, $\nu_3$, $\nu_4$. Assume that end points of the arcs $\nu_1$, $\nu_2$, $\nu_3$, $\nu_4$ on $N^-$ are lined up consecutively as $N^-\cap\nu_1$, $N^-\cap\nu_2$, $N^-\cap\nu_3$, $N^-\cap\nu_4$. The curve $c_S$ has two arcs $a_1$, $a_2$ on $T^-$ of slope 0 and  two arcs $b_1$, $b_2$ on $T^+$ of slope $\infty$. Assume that $\nu_i\cap a_1\neq \emptyset$ for $i=1,2$ and  $\nu_1\cap b_1\neq \emptyset$. See figure \ref{fig:IA2b}. There are eight regions $D_1$,...,$D_8$ on $N$ that can contain a vertex of a bigon.  The regions  $D_1$,...,$D_8$ are shown in figure \ref{fig:IA2b}. 
\begin{figure}[h]
\begin{center}
\psfrag{d1}[][][0.7]{$D_1$}
\psfrag{d2}[][][0.7]{$D_2$}
\psfrag{d3}[][][0.7]{$D_3$}
\psfrag{d4}[][][0.7]{$D_4$}
\psfrag{d5}[][][0.7]{$D_5$}
\psfrag{d6}[][][0.7]{$D_6$}
\psfrag{d7}[][][0.7]{$D_7$}
\psfrag{d8}[][][0.7]{$D_8$}
\psfrag{a1}[][][1]{$a_1$}
\psfrag{a2}[][][1]{$a_2$}
\psfrag{b1}{$b_1$}
\psfrag{b2}{$b_2$}
\psfrag{n-}[][][1]{$N^-$}
\psfrag{n+}[][][1]{$N^+$}
\psfrag{n1}[][][1]{$\nu_1$}
\psfrag{n2}[][][1]{$\nu_2$}
\psfrag{n3}[][][1]{$\nu_3$}
\psfrag{n4}[][][1]{$\nu_4$}
\psfrag{cp}[][][1]{$c_P$}
\includegraphics[height=5.5cm, width=5.5cm]{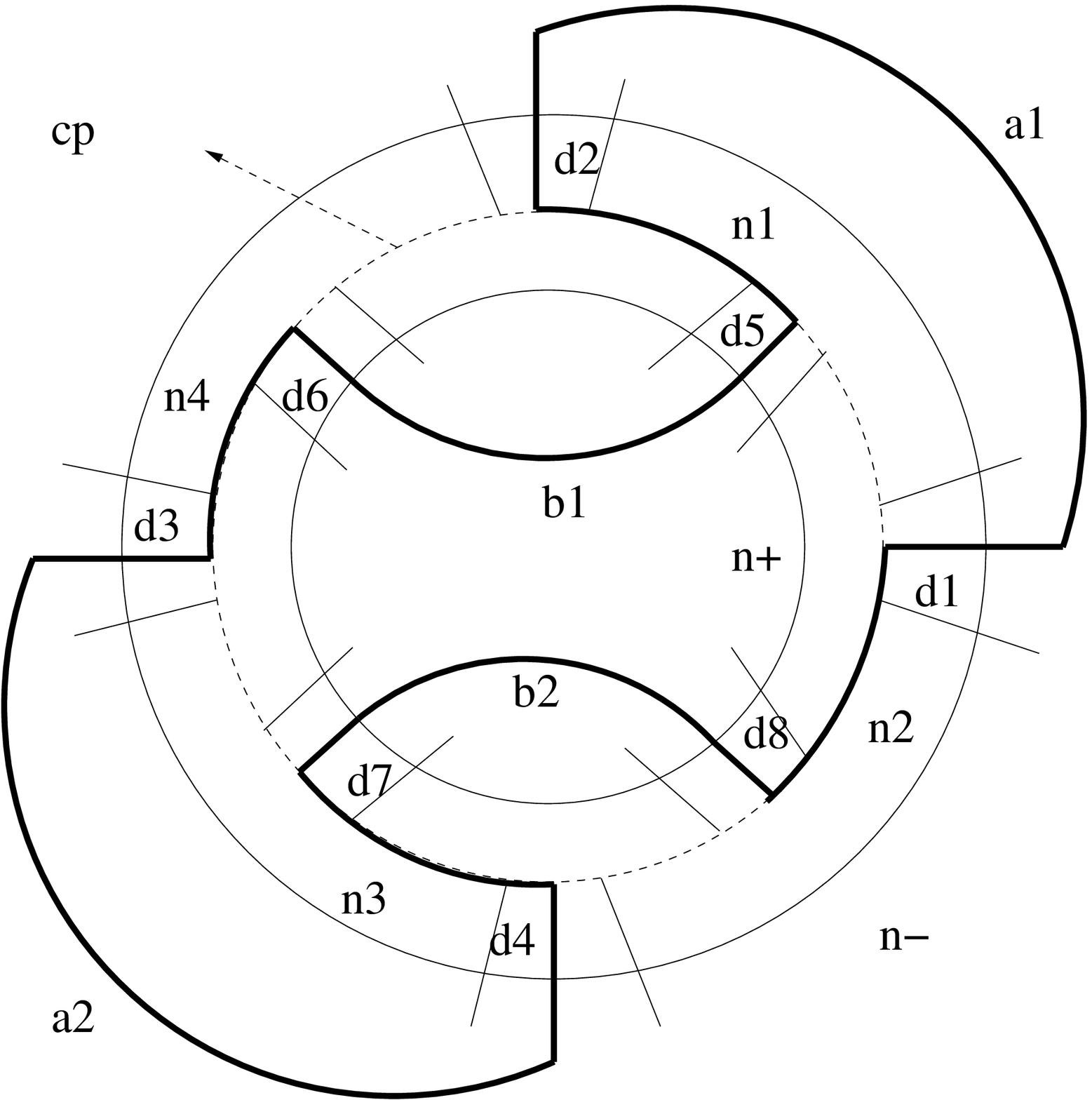}
\end{center}
\caption{}\label{fig:IA2b}
\end{figure}
Any bigon should contain two of them.  After an isotopy, we may assume that $\alpha(c_Q)=c_Q$ and  $\alpha(c_S)=c_S$. Then  $\alpha(D_i)=D_{i+2}$ for $i=1,2$ and $\Theta(\{D_i\ |\ i=1,...,4\})=\{D_i\ |\ i=5,...,8\}$ (see Definition \ref{defn:Theta} for $\Theta$). So it is enough to check if $D_i$ is a part of a bigon for $i=1,2$.

\begin{description}
\item[$D_1$] The region $D_1$ is part of a region $\widetilde{D}_1$ in $T$ whose four consecutive sides are $x_1$, $a_1$, $x_2$, $x_3$ where $x_1\in F^+_{Q,k}\cup F^+_{Q,k+1}$, $x_2\in F^-_{Q,0}$ and $x_3\in F^+_{Q,k}\cup F^+_{Q,k+1}$. See figure \ref{fig:IA2d}(a). If  $\widetilde{D}_1$ is a bigon then $v_Q\cdot v_P<2(n+m)$, a contradiction. 
\item[$D_2$]
\begin{itemize}
\item If $b=0$ then $a\neq 0$.  Then  $D_2$ is part of a region $\widetilde{D}_2$ whose eight sides are  $x_1$, $y_1$, $a_1$, $z_1$, $x_2$, $a_2$, $y_2$, $z_2$ where $x_1,x_2,z_1,z_2\in F^+_{Q,\infty}$, \  $y_1,y_2\in F^-_{Q,1/k}$. See figure \ref{fig:IA2d}(b). Therefore $\tilde{D}_1$ can not be a bigon. 

\item If $a,b\neq 0$ then $D_2$ is part of a region $\widetilde{D}_2$ whose four sides  are $x_1$, $a_2$, $y_1$, $z_1$ where $x_1\in F^+_{Q,\infty}$, $y_1\in F^-_{Q,1/(k+1)}$, and $z_1$ is either a piece of  $c_S$ or $z_1\in F^+_{Q,k}\cup F^+_{Q,k+1}$. See figure \ref{fig:IA2d}(c).  If $z_1$ is an arc in $c_S$ then $\widetilde{D}_2$ can not be a bigon. If  $z_1\in F^+_{Q,\infty}\cup F^+_{Q,k}\cup F^+_{Q,k+1}$ and  $\widetilde{D}_2$ is a bigon then $v_Q\cdot v_P<2(n+m)$, a contradiction. 
\end{itemize}

\begin{figure}[h]
\begin{center}
\psfrag{d1}[][][1]{$D_1$}
\psfrag{d2}[][][1]{$D_2$}
\psfrag{a1}[][][1]{$a_1$}
\psfrag{a2}[][][1]{$a_2$}
\psfrag{z1}[][][1]{$z_1$}
\psfrag{z2}[][][1]{$z_2$}
\psfrag{y1}[][][1]{$y_1$}
\psfrag{y2}[][][1]{$y_2$}
\psfrag{x1}[][][1]{$x_1$}
\psfrag{x2}[][][1]{$x_2$}
\psfrag{x3}[][][1]{$x_3$}
\psfrag{(a)}[][][1]{$(a)$}
\psfrag{(b)}[][][1]{$(b)$}
\psfrag{(c)}[][][1]{$(c)$}
\includegraphics[height=3.5cm, width=3cm]{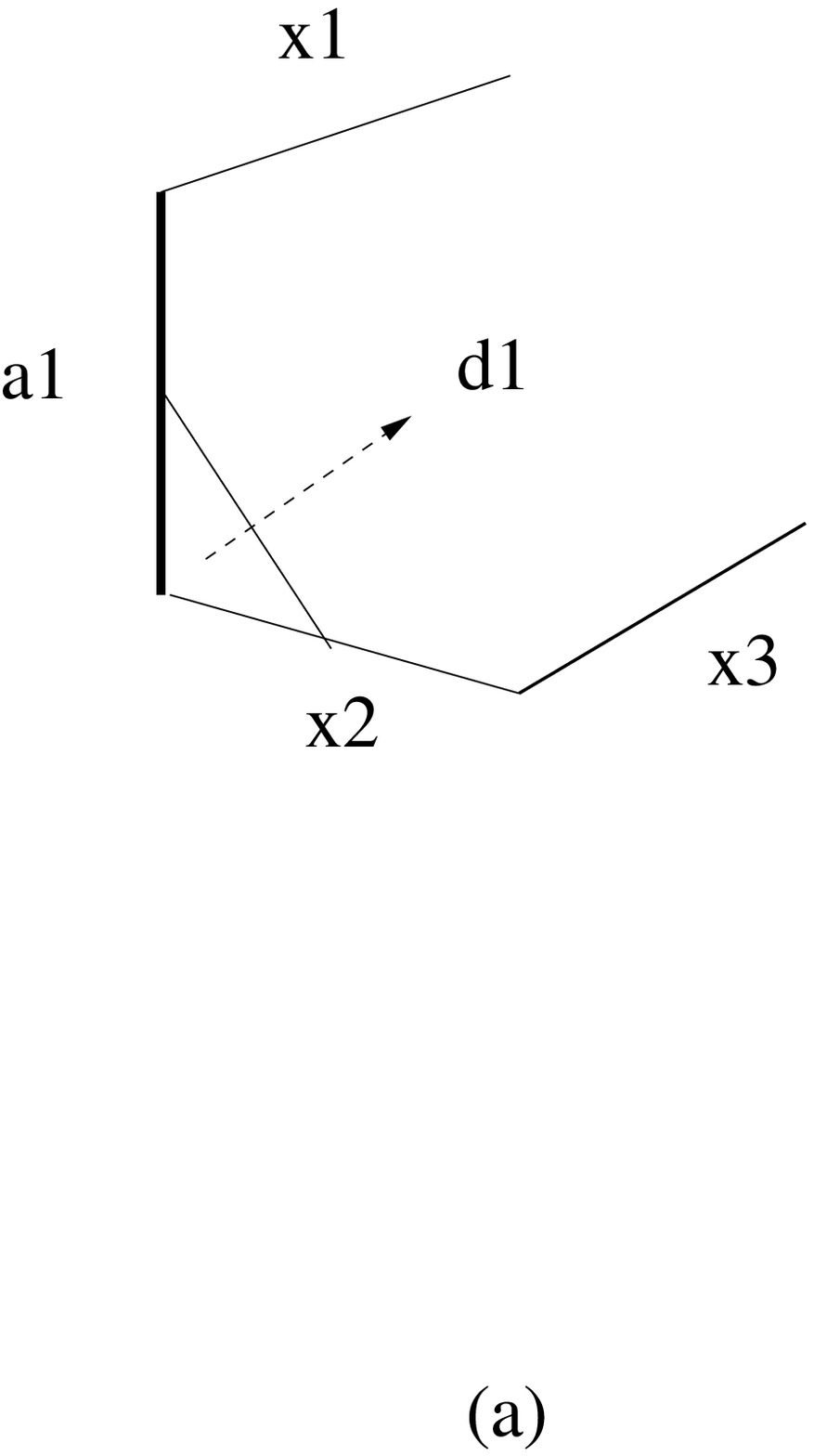}\hskip2cm
\includegraphics[height=3.5cm, width=4cm]{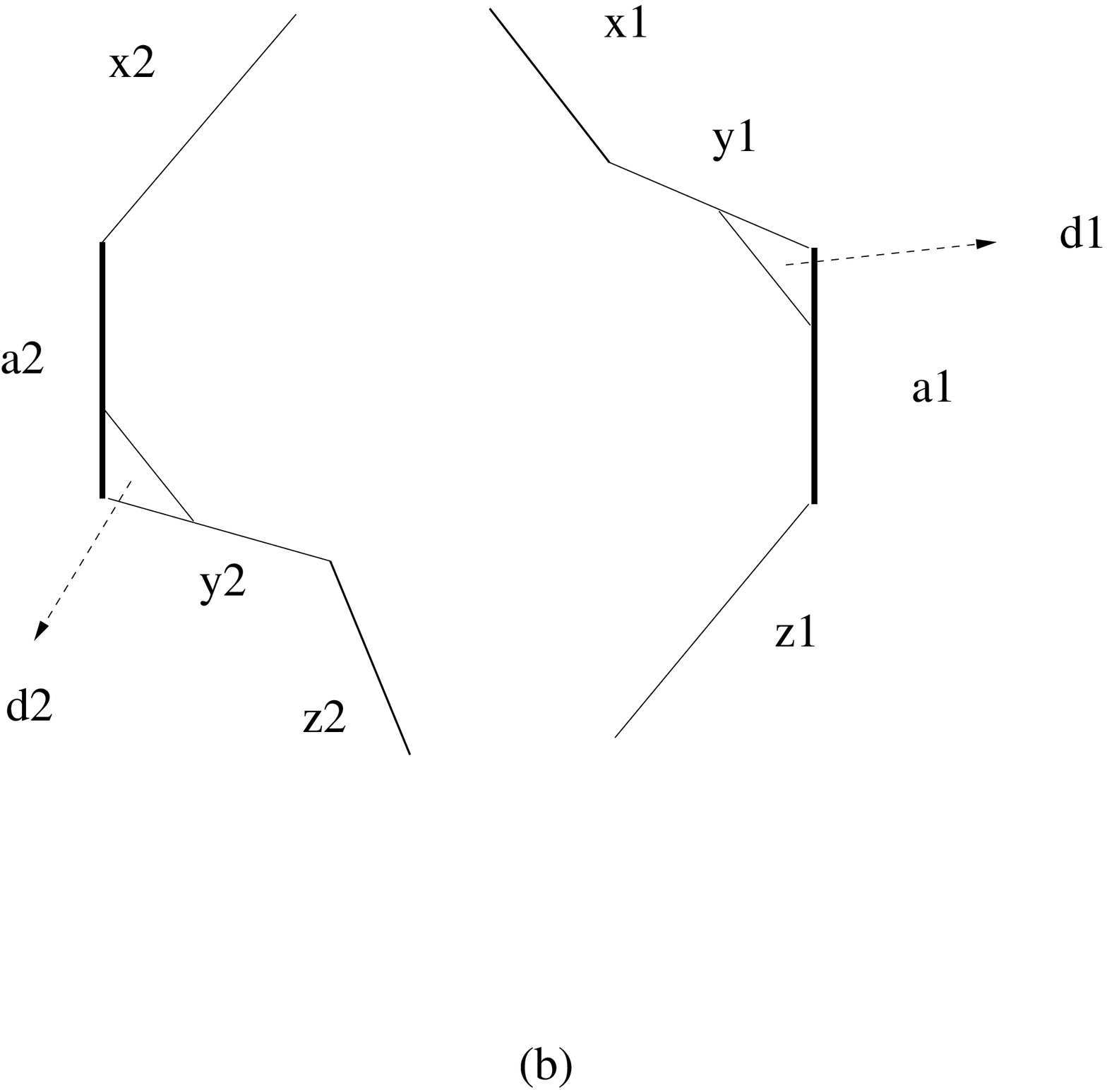}\hskip2cm
\includegraphics[height=3.5cm, width=3cm]{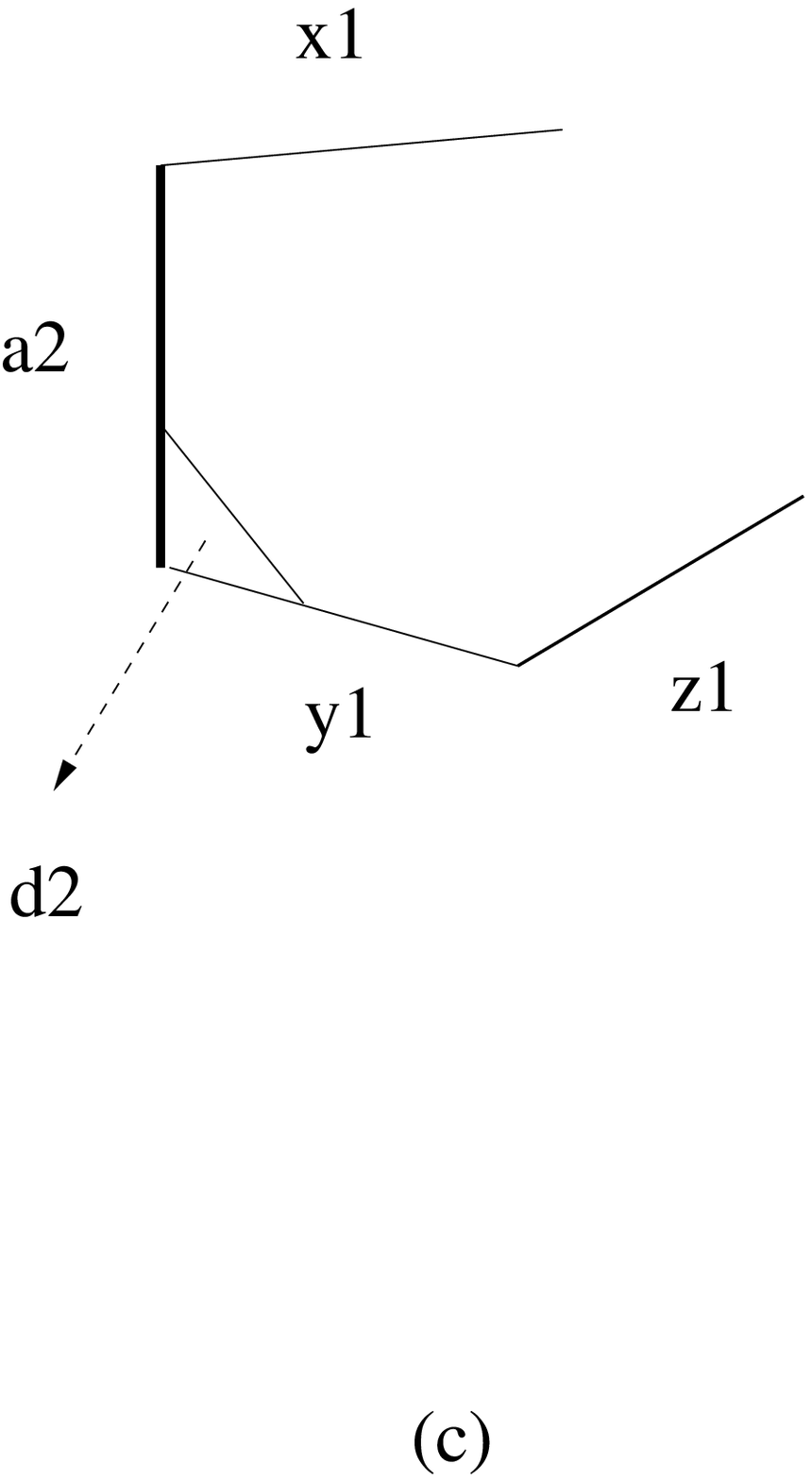}
\end{center}
\caption{}\label{fig:IA2d}
\end{figure}
\end{description}

By the above cases, $v_S\cdot v_Q=|c_S\cap c_Q|$.

In figure \ref{fig:IA2e},  intersection of a reducing sphere $R^\prime$ with the surface $T$ is shown. 
Notice that $R^\prime\in v_{\gamma S}$ and $v_{S}\cdot v_{\gamma S}=4$.  By an argument similar to the proof of Claim 1 we can show that $v_{R^\prime}\cdot v_Q=|c_{R^\prime}\cap c_Q|=4kb+4(k-1)a+2m+2n=v_{\gamma S} \cdot v_Q\geq2m+2n$. 
\begin{figure}[h]
\begin{center}
\psfrag{rp}[][][1.5]{$c_{R^\prime}$}
\psfrag{f-0}[][][1]{$F^-_{Q,0}$}
\psfrag{f-1k}[][][1]{$F^-_{Q,\frac{1}{k}}$}
\psfrag{f-1k+1}[][][1]{$F^-_{Q,\frac{1}{k+1}}$}
\psfrag{f+i}[][][1]{$F^+_{Q,\infty}$}
\psfrag{f+k}[][][1]{$F^+_{Q,k}$}
\psfrag{f+k+1}[][][1]{$F^+_{Q,k+1}$}
\includegraphics[height=6cm, width=6cm]{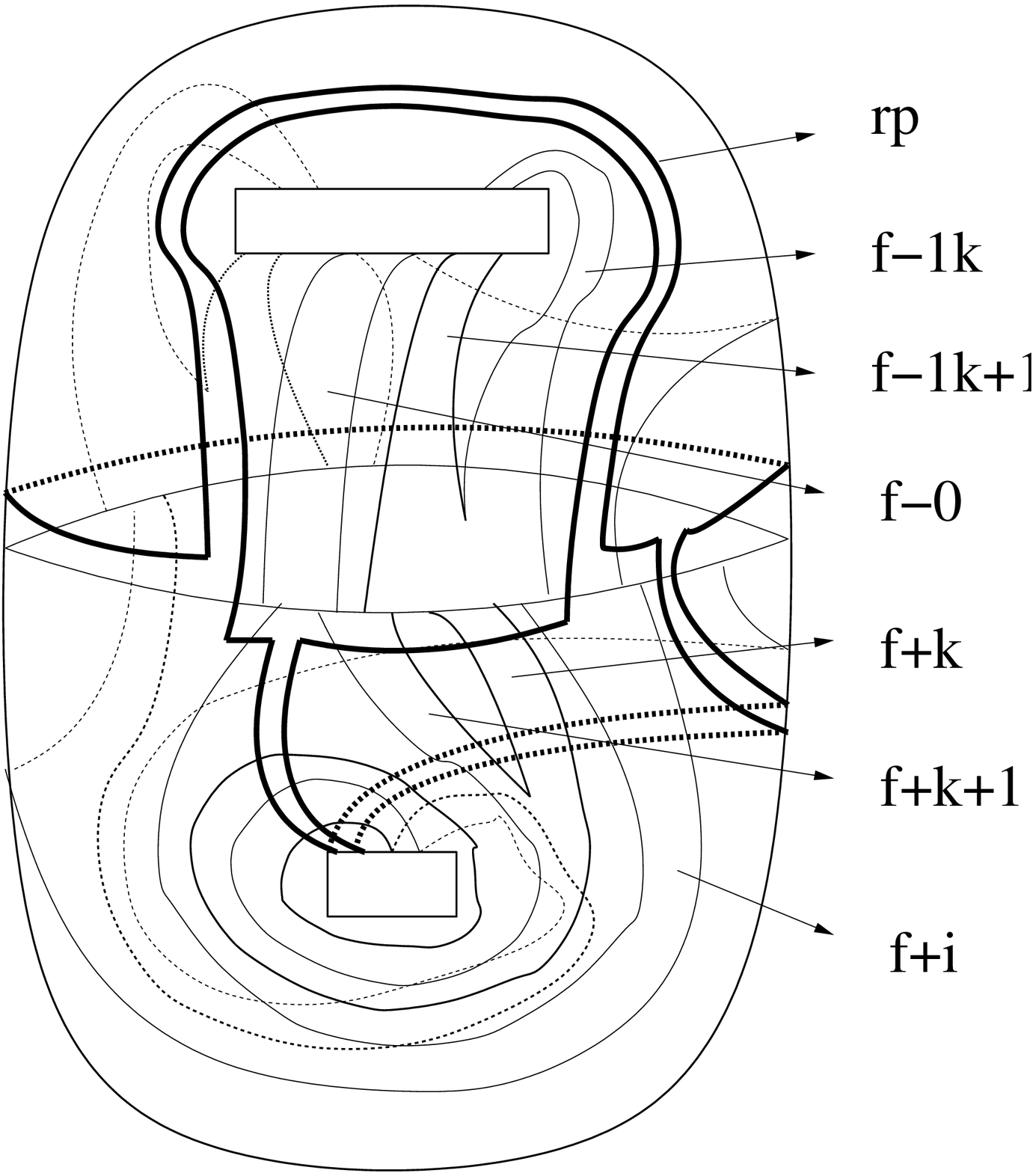}
\end{center}
\caption{}\label{fig:IA2e}
\end{figure}

\textbf{Proof of  Claim 2.} We will do the calculation for $i=\pm1$. The general case is similar.
We may assume that  $\beta^i (c_S)$ and $\beta^i\gamma(c_S)$ intersect  $c_Q$ in a neighborhood $N$ described in the proof of Claim 1. By an argument similar to the proof of Claim 1 we get
\begin{itemize}
\item  $v_{\beta(S)} \cdot v_Q=4p+2m+6n> 2(n+m)$. See figure \ref{fig:IA2f} (a).
\item $v_{\beta^{-1}(S)} \cdot v_Q=6m+2n-4p> 2(n+m)$.  See figure \ref{fig:IA2f} (b).
\begin{figure}[h]
\begin{center}
\psfrag{b+r}[][][1]{$\beta (c_S)$}
\psfrag{n}[][][1]{$N$}
\psfrag{bcs}[][][1]{$\beta(c_S)$}
\psfrag{b-r}[][][1]{$\beta^{-1}(c_S)$}
\psfrag{(a)}[][][1]{$(a)$}
\psfrag{(b)}[][][1]{$(b)$}
\includegraphics[height=5.5cm, width=4.5cm]{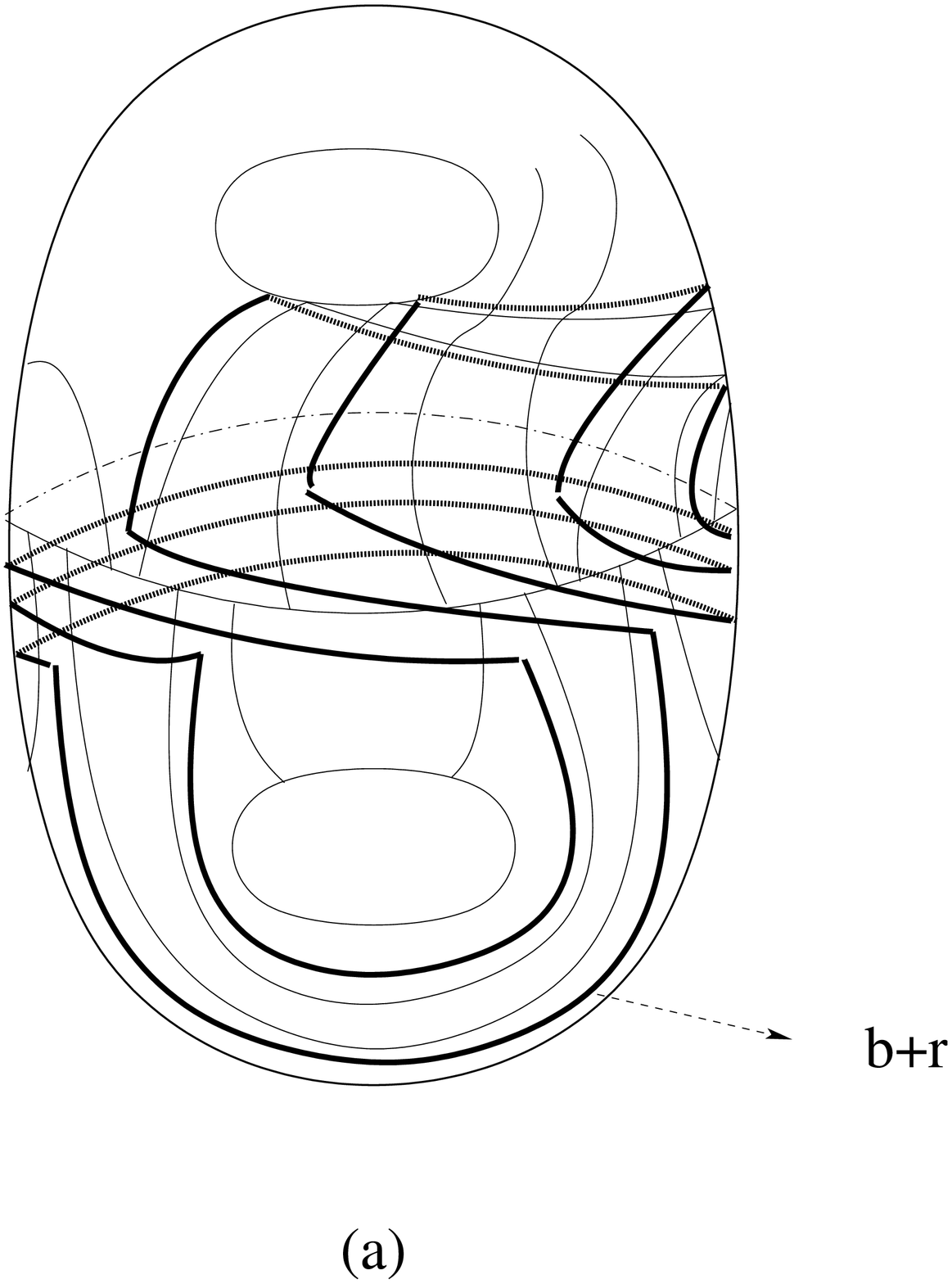}\hskip1cm
\includegraphics[height=5.5cm, width=4.5cm]{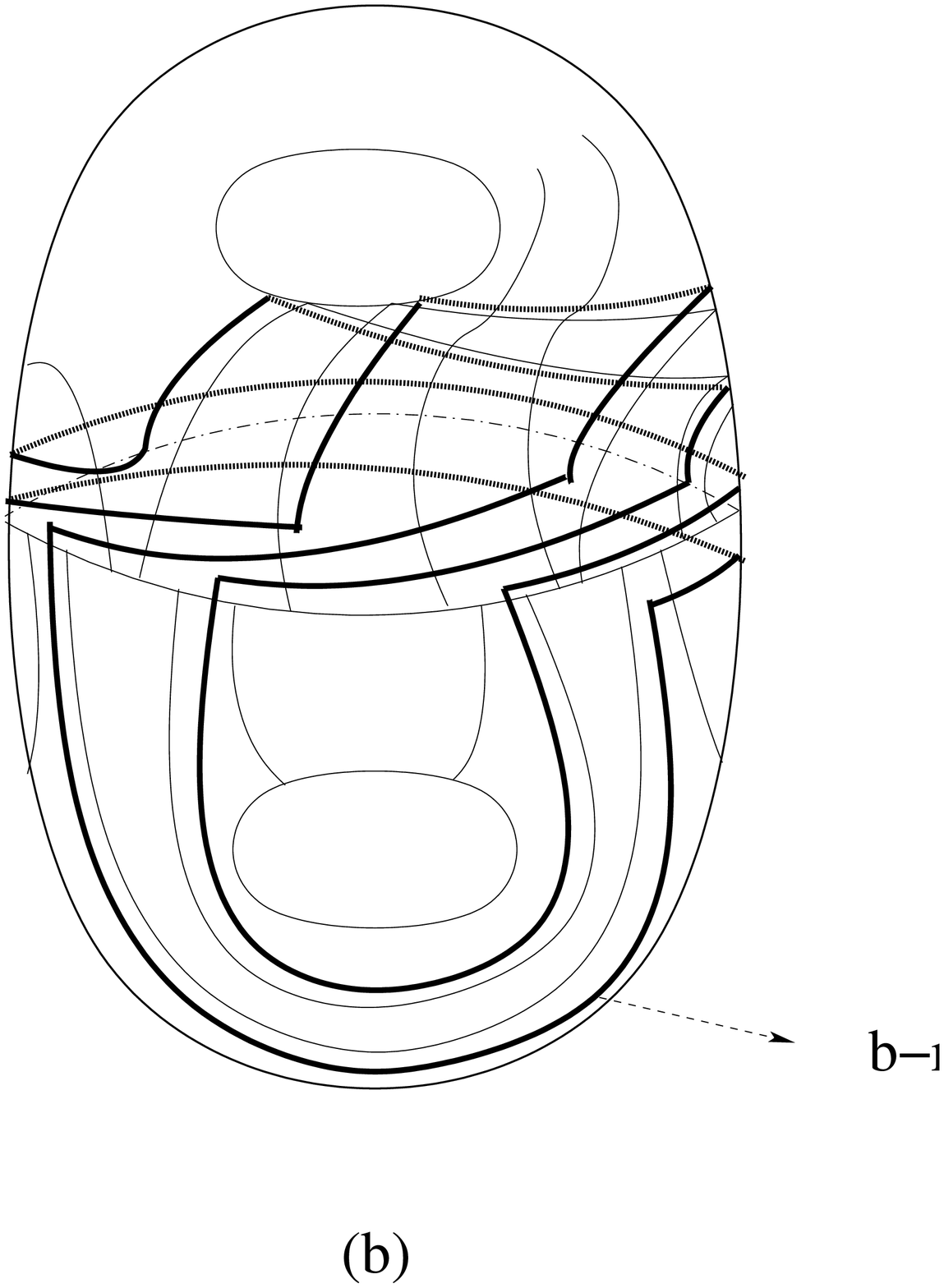}
\end{center}
\caption{}\label{fig:IA2f}
\end{figure}
\item $v_{\beta\gamma(S)}\cdot v_Q =4kb+4(k-1)a+4m+2n+2p> 2(n+m)$. See figure \ref{fig:IA2g} (a).
\item $v_{\beta^{-1}\gamma(S)}\cdot v_Q =4kb+4(k-1)a+6m+6n-4p> 2(n+m)$. See figure \ref{fig:IA2g} (b).
\begin{figure}[h]
\begin{center}
\psfrag{b+rp}[][][1]{$\beta(c_{R^\prime})$}
\psfrag{b-rp}[][][1]{$\beta^{-1}(c_{R^\prime})$}
\psfrag{(a)}[][][1]{$(a)$}
\psfrag{(b)}[][][1]{$(b)$}
\includegraphics[height=6cm, width=4.5cm]{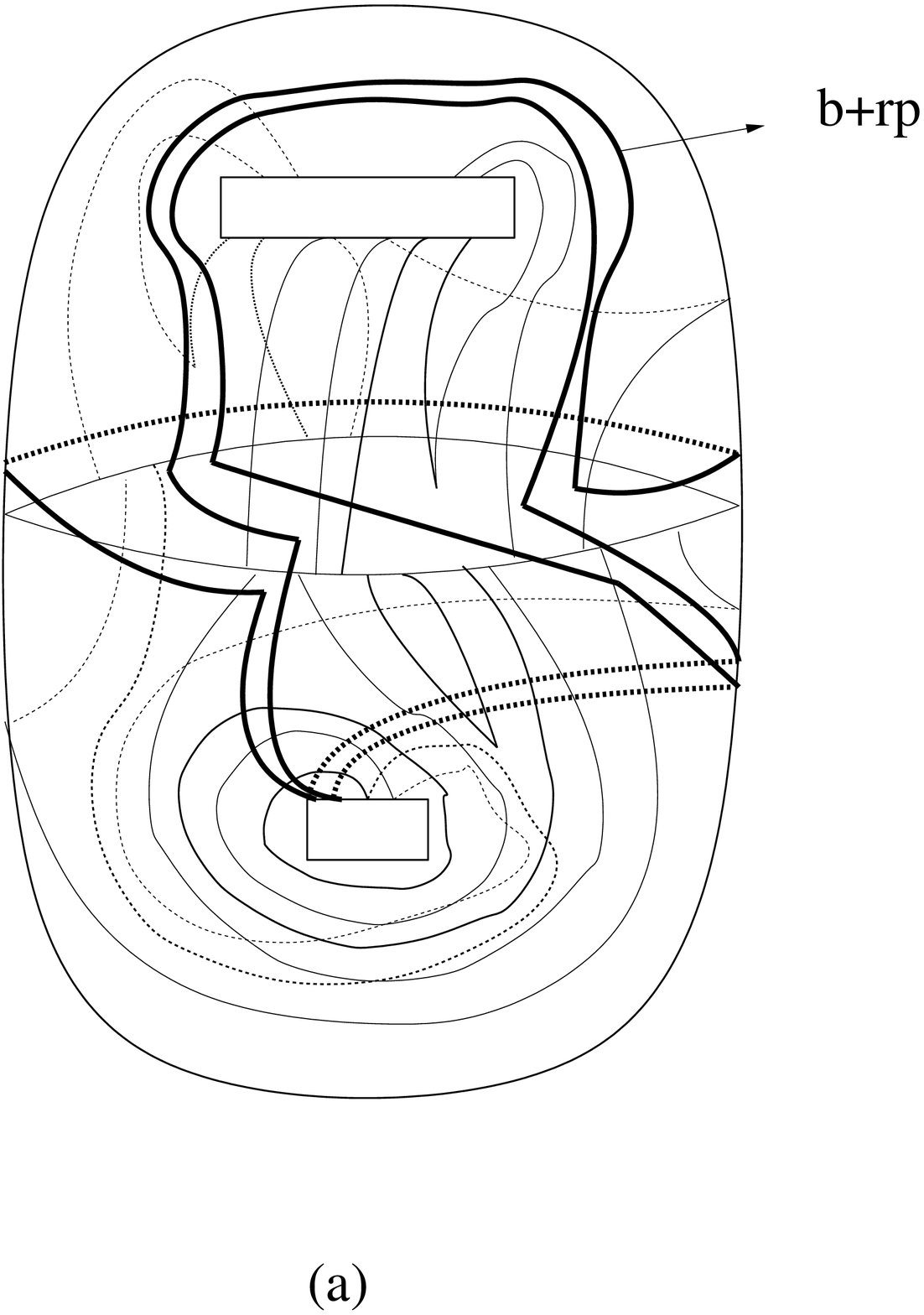}\ \ \ \ \ \ \ 
\includegraphics[height=6cm, width=4.5cm]{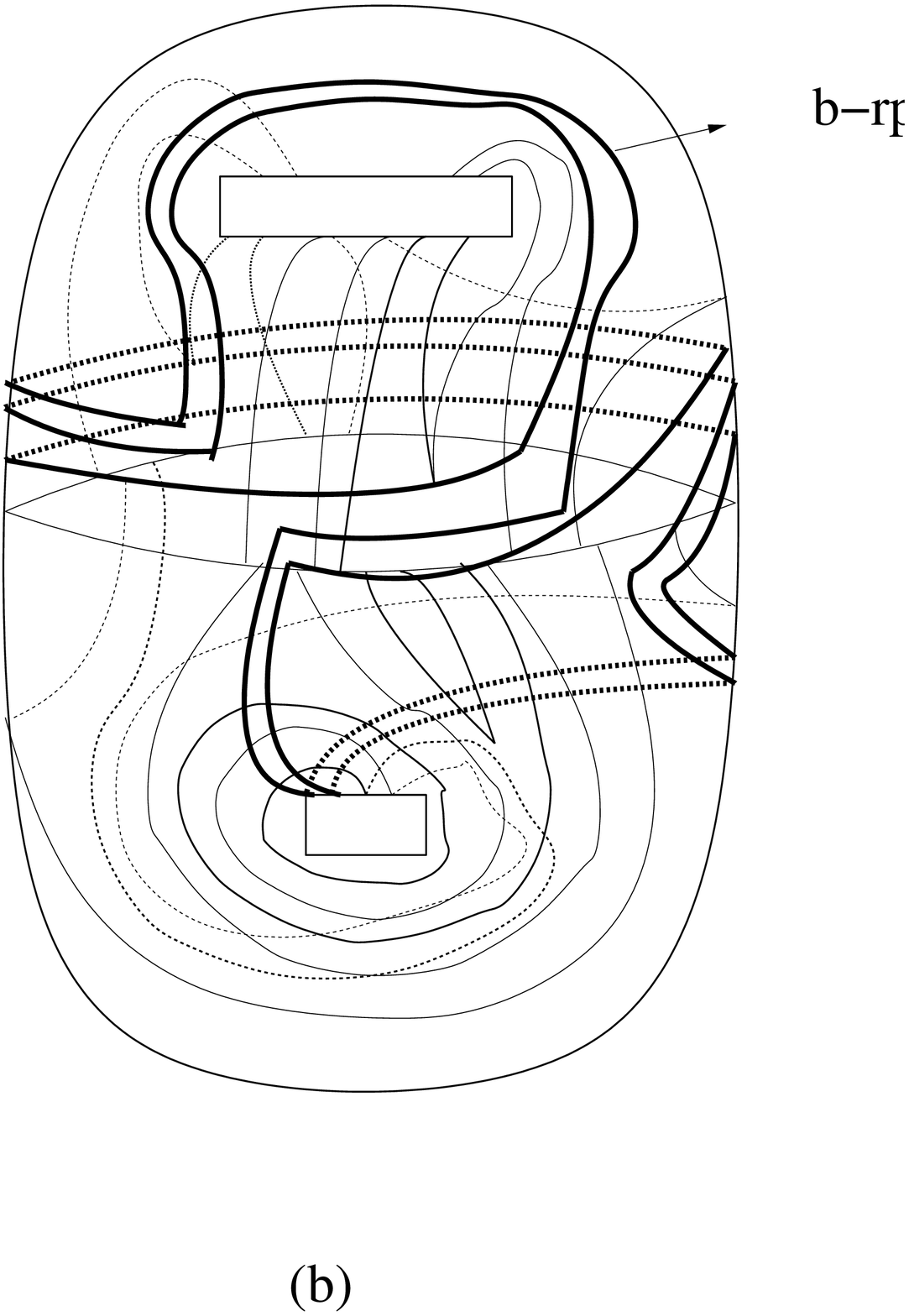}
\end{center}
\caption{}\label{fig:IA2g}
\end{figure}
\end{itemize}

This implies that the vertex $v_R=v_S$ and satisfies the conditions  of Proposition \ref{backbone}.

\newpage
\noindent {\bf I.A.2.}  If $\{d_{i}\}\subseteq \{e_{ij}\}$ (see  figure \ref{fig:IA3a}):  Set $p=|\{e_{0j}\}\cap \{h_{0j}\}|$. Then $0<p\leq m-n$. Either $p<m-n-p$ or $m-n-p<p$. Assume $p<m-n-p$. Consider the curve $\xi$ shown in figure \ref{fig:IA3a}. The curve $\xi$ is an intersection of a reducing sphere $S$ with $T$. Denote $\xi$ by $c_S$. Notice that $v_S\cdot v_P=4$. 
\begin{figure}[h]
\begin{center}
\psfrag{fi-}[][][1]{$F^-_{Q,\infty}$}
\psfrag{f1-}[][][1]{$F^-_{Q,1}$}
\psfrag{f0-}[][][1]{$F^-_{Q,0}$}
\psfrag{f1k-}[][][1]{$F^-_{Q,\frac{1}{k}}$}
\psfrag{f1k+1-}[][][1]{$F^-_{Q,\frac{1}{k+1}}$}
\psfrag{fi+}[][][1]{$F^+_{Q,\infty}$}
\psfrag{f1+}[][][1]{$F^+_{Q,1}$}
\psfrag{f0+}[][][1]{$F^+_{Q,0}$}
\psfrag{fk1+}[][][1]{$F^+_{Q,k}$}
\psfrag{fk+11+}[][][1]{$F^+_{Q,k+1}$}
\psfrag{cr}[][][1.75]{$\xi$}
\psfrag{mix1}{$F^-_{Q,\frac{1}{k}}\cup F^-_{Q,\frac{1}{k+1}}$}
\psfrag{mix2}{$F^+_{Q,k}\cup F^+_{Q,k+1}$}
\includegraphics[height=5cm, width=6cm]{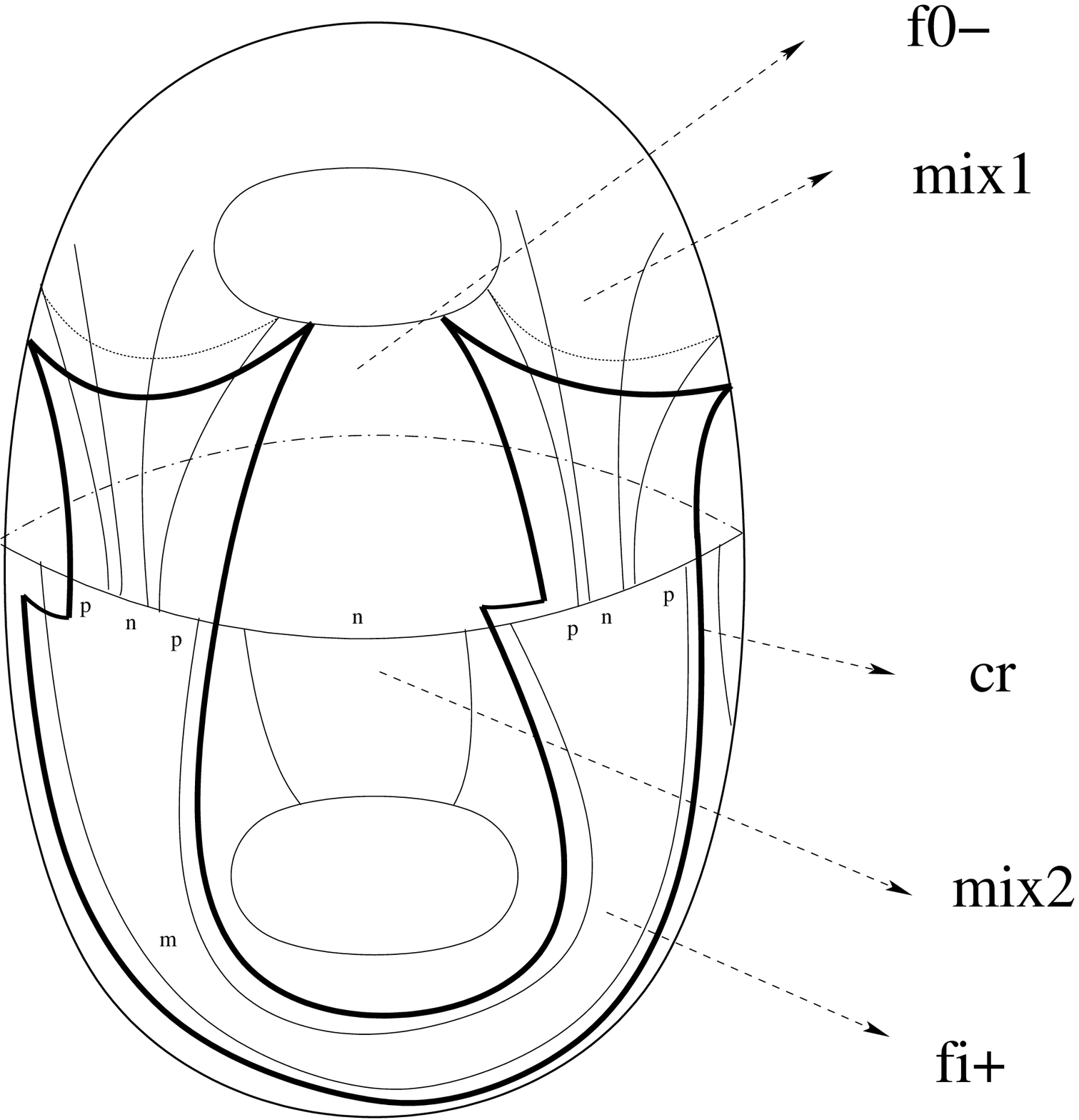}
\end{center}
\caption{}\label{fig:IA3a}
\end{figure}

By an argument similar to the proof of the case I.A.1. we get 
\begin{itemize}
\item $ v_S\cdot v_Q=|c_S\cap c_Q|=2(m-n-2p) < v_P\cdot v_Q=2(n+m)$
\item $v_S \cdot v_{\gamma(S)} =4$
\item $v_{\gamma(S)}\cdot v_Q =4kb+4(k-1)a+2(m+n)\geq 2(m+n)$ (see figure \ref{fig:IA3b})
\item $v_{\beta^i (S)} \cdot v_Q,\ \ v_{\beta^i\gamma(S)}\cdot v_Q > 2(n+m)$ for $i\neq 0$.
\end{itemize}
This implies that the vertex $v_R=v_S$ and satisfies the conditions  of Proposition \ref{backbone}.
\begin{figure}[h]
\begin{center}
\psfrag{fi-}[][][1]{$F^-_{Q,\infty}$}
\psfrag{f1-}[][][1]{$F^-_{Q,1}$}
\psfrag{f0-}[][][1]{$F^-_{Q,0}$}
\psfrag{f1k-}[][][1]{$F^-_{Q,\frac{1}{k}}$}
\psfrag{f1k+1-}[][][1]{$F^-_{Q,\frac{1}{k+1}}$}
\psfrag{fi+}[][][1]{$F^+_{Q,\infty}$}
\psfrag{f1+}[][][1]{$F^+_{Q,1}$}
\psfrag{f0+}[][][1]{$F^+_{Q,0}$}
\psfrag{fk1+}[][][1]{$F^+_{Q,k}$}
\psfrag{fk+11+}[][][1]{$F^+_{Q,k+1}$}
\psfrag{crp}[][][1.75]{$c_{R^\prime}$}
\psfrag{mix1}{$F^-_{Q,\frac{1}{k}}\cup F^-_{Q,\frac{1}{k+1}}$}
\psfrag{mix2}{$F^+_{Q,k}\cup F^+_{Q,k+1}$}
\includegraphics[height=5.5cm, width=6cm]{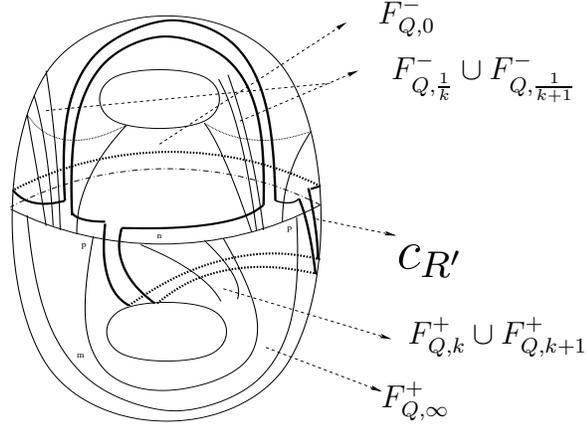}
\end{center}
\caption{The curve $c_{R^\prime}$ in the figure is $R^\prime\cap T$ for some reducing sphere $R^\prime$ for $T$ satisfying $R^\prime\in v_{\gamma S}$}\label{fig:IA3b}
\end{figure}

\noindent {\bf I.B.}  If $n=0$: This is a special case of I.A.3.\\

\noindent{\bf Case II.} If  $N(Q,T^-,0)=m$, $N(Q,T^-,\infty)=n\neq 0=N(Q,T^-,1)$ then $N(Q,T^+,0)=n$, $N(Q,T^+,\infty)=m\neq 0=N(Q,T^+,1)$ by Proposition \ref{a=1/a}. By Lemma \ref{notequal}, $m\neq n$. Suppose $m<n$.   By Lemma \ref{lemmaC},  $ \{ e_{ij} |i=0,1 \ j=1,..,m \} \subseteq \{f_{ij} |i=0,1 \ j=2,..,n-1\}$. By the argument in \cite[Lemma 5]{Sc}  we get  two non-isotopic reducing  spheres for $T$ that satisfy  (i) and (ii).  Let us call $S$ the one having an arc on  $T^-$  of slope $0$ and $S^\prime$ the one having an arc on  $T^+$ of slope $0$.  In the figure \ref{fig:II} intersections of two reducing spheres $R$ and $R^\prime$ with $T$ are shown. It is easy to see that $R\in v_S$ and  $R^\prime\in v_{S^\prime}$.
\begin{figure}[h]
\begin{center}
\psfrag{cr}[][][1.75]{$c_R$}
\psfrag{crp}[][][1.75]{$c_{R^\prime}$}
\includegraphics[height=6cm, width=6cm]{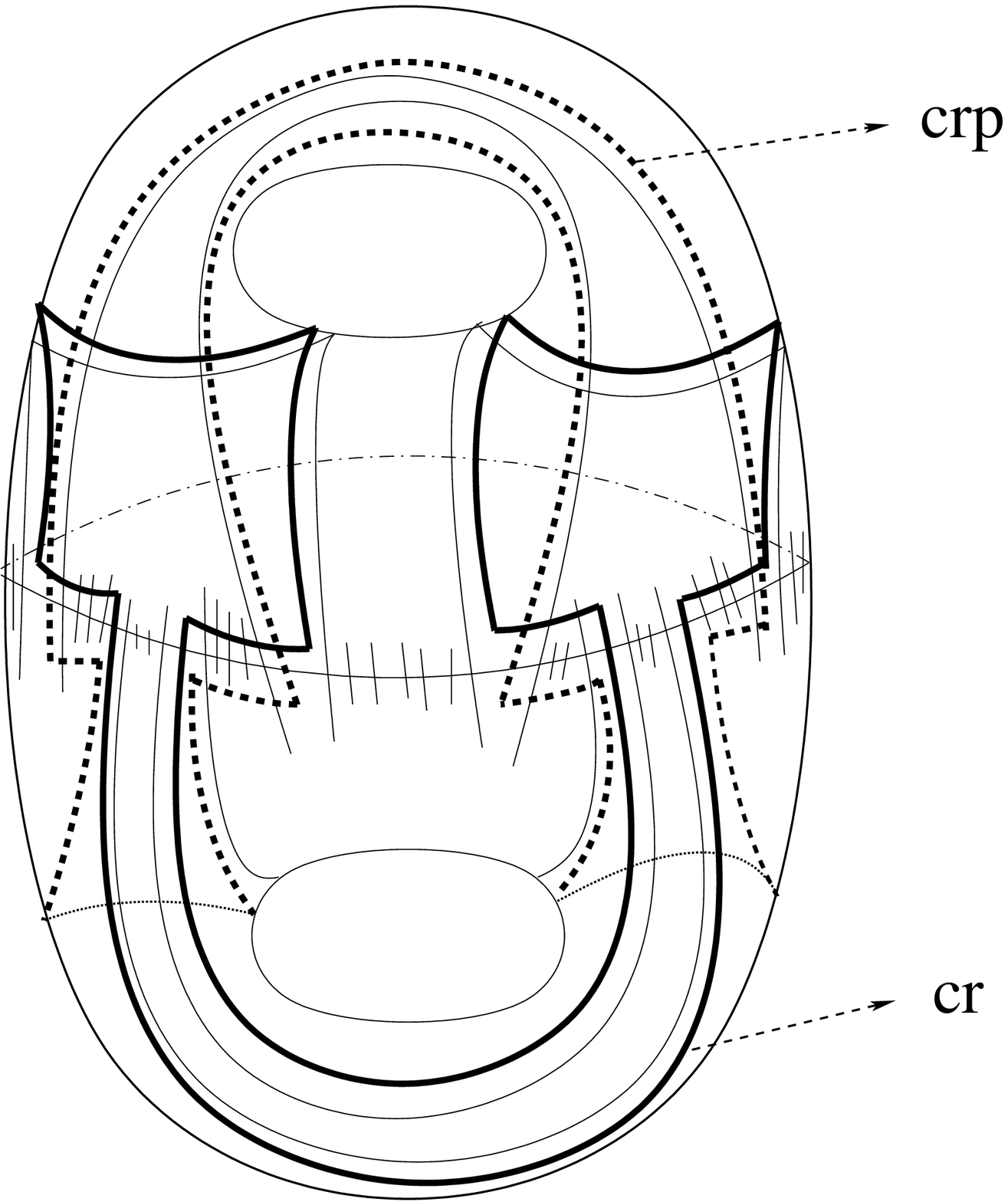}
\end{center}
\caption{}\label{fig:II}
\end{figure}

Let $p=|\{g_{0j}\}\cap \{f_{0j}\}|$. Then $0<p\leq m-n$. Either $p<m-n-p$ or $m-n-p<p$. Assume $p<m-n-p$. Then  by an argument similar to the proof of the case I.A.1. we can show that $2n+2m = v_P\cdot v_Q > v_R\cdot v_Q = 2n-2m > v_{R^\prime} \cdot v_Q =2(n-m-2p)$, $v_{R} \cdot v_{R^\prime}=4$ and $v_{\beta^i (R)} \cdot v_Q, \  v_{\beta^i(R^\prime)}\cdot v_Q > 2n+2m$ for $i\neq 0$.

{\bf Case III.} If  $N(Q,T^-,0)=m$, $N(Q,T^-,\infty)=n$,  $N(Q,T^-,1)=p$ where $m,n,p\neq 0$,  then $N(Q,T^+,0)=n$, $N(Q,T^+,\infty)=m$, $N(Q,T^+,1)=p$ by Proposition \ref{a=1/a}. By Lemma \ref{notequal}, $m\neq n$. Say $m>n$. 

The curves $A$, $B$, $C$ and  $c_P$ divide $T$ into four punctured discs  $T^-_f$, $T^-_b$, $T^+_f$, $T^+_b$ where  $T^-_f\cup T^-_b= T^-$  and $T^+_f\cup T^+_b= T^+$. This division also gives two pairs of pants $T^-_f\cup T^+_f=P_f$ and $T^-_b\cup T^+_b= P_b$. Let $c_f= P_f\cap c_P$ and $c_b= P_b\cap c_P$. 

Let $K$ be a reducing sphere intersecting  the interior of $T^-$  in a simple  closed curve parallel to  $c_P$.  The reducing sphere $K$ divides $T$  into two parts.  Denote the one  containing the curve $B$  by $t^-$ and the one containing the curve $C$ by $t^+$.  Let  $c^f_{K}= T^-_f\cap K$ and $c^b_{K}= T^-_b\cap K$.

Suppose that
$F^-_{Q,0} \cap t^- \cap A = F^-_{Q,1} \cap t^- \cap A = \emptyset \ $,
$|F^-_{Q,\infty}\cap (c^f_{K}\setminus A)|$ $=$
$|F^-_{Q,\infty}\cap (c^b_{K}\setminus A)|$ $=$
$|F^-_{Q,\infty}\cap t^-\cap A|=n$ and that
$k^\prime_{01}$, $k^\prime_{02}$,..., $k^\prime_{0p}$,
$e^\prime_{01}$, $e^\prime_{02}$,..., $e^\prime_{0m}$ and
$g^\prime_{01}$, $g^\prime_{02}$,..., $g^\prime_{0n}$ are consecutive intersection points of the arcs in $F^-_{Q,1}$, $F^-_{Q,0}$ and $F^-_{Q,\infty}$ with $c^f_{K}$ respectively. Locate arcs of $c_Q$ on $T^+$ in such a way that
$|F^+_{Q,\infty}\cap (c^f_{P}\setminus A)|$ $=$
$|F^+_{Q,\infty}\cap (c^b_{P}\setminus A)|$ $=$
$|F^+_{Q,\infty}\cap A|=m$ and
$|F^+_{Q,0}\cap A|=|F^+_{Q,1}\cap A|=0$. Suppose that $l_{01}$,..., $l_{0p}$, $f_{01}$,..., $f_{0n}$ and $h_{01}$,..., $h_{0m}$ are consecutive intersection points of the arcs in $F^+_{Q,1}$, $F^+_{Q,0}$ and $F^+_{Q,\infty}$  with  $c_f$ respectively. Let $\tau$ be an arc in $F^-_{Q,1}$ whose intersection with $c^f_{K}$ is $k^\prime_{01}$.  Suppose that $\tau \cap (t^+\setminus T^+) \cap A\neq \emptyset$.   See figure \ref{fig:III}.  By applying a power of $\beta$ we can assume that  $2 \leq |c_Q \cap A \cap (t^+\setminus T^+)|<2(p+n+m) $.  By the argument in \cite[Lemma 5]{Sc}  we get  two non-isotopic reducing  spheres for $T$ that satisfy  (i) and (ii).  Let us call $S$ the one having an arc on  $T^-$  of slope $0$ and $S^\prime$ the one having an arc on  $T^+$ of slope $0$. 
\begin{figure}[h]
\begin{center}
\psfrag{fi-f}[][][1]{$F^-_{Q,\infty} \cap t^- \cap P_f$}
\psfrag{f1-f}[][][1]{$F^-_{Q,1}\cap t^- \cap P_f$}
\psfrag{f0-f}[][][1]{$F^-_{Q,0}\cap t^- \cap P_f$}
\psfrag{fi+f}[][][1]{$F^+_{Q,\infty}\cap t^+ \cap P_f$}
\psfrag{f1+f}[][][1]{$F^+_{Q,1}\cap t^+ \cap P_f$}
\psfrag{f0+f}[][][1]{$F^+_{Q,0}\cap t^+ \cap P_f$}
\psfrag{3}[][][1]{$P_f$}
\psfrag{cf}[][][1]{$c_f$}
\psfrag{cfpp}[][][1]{$c^f_{K}$}
\psfrag{tau}[][][1]{$\tau$}
\psfrag{B}[][][1]{$B$}
\psfrag{C}[][][1]{$C$}
\psfrag{m}[][][1]{$m$}
\psfrag{n}[][][1]{$n$}
\psfrag{p}[][][1]{$p$}
\includegraphics[height=7cm, width=10cm]{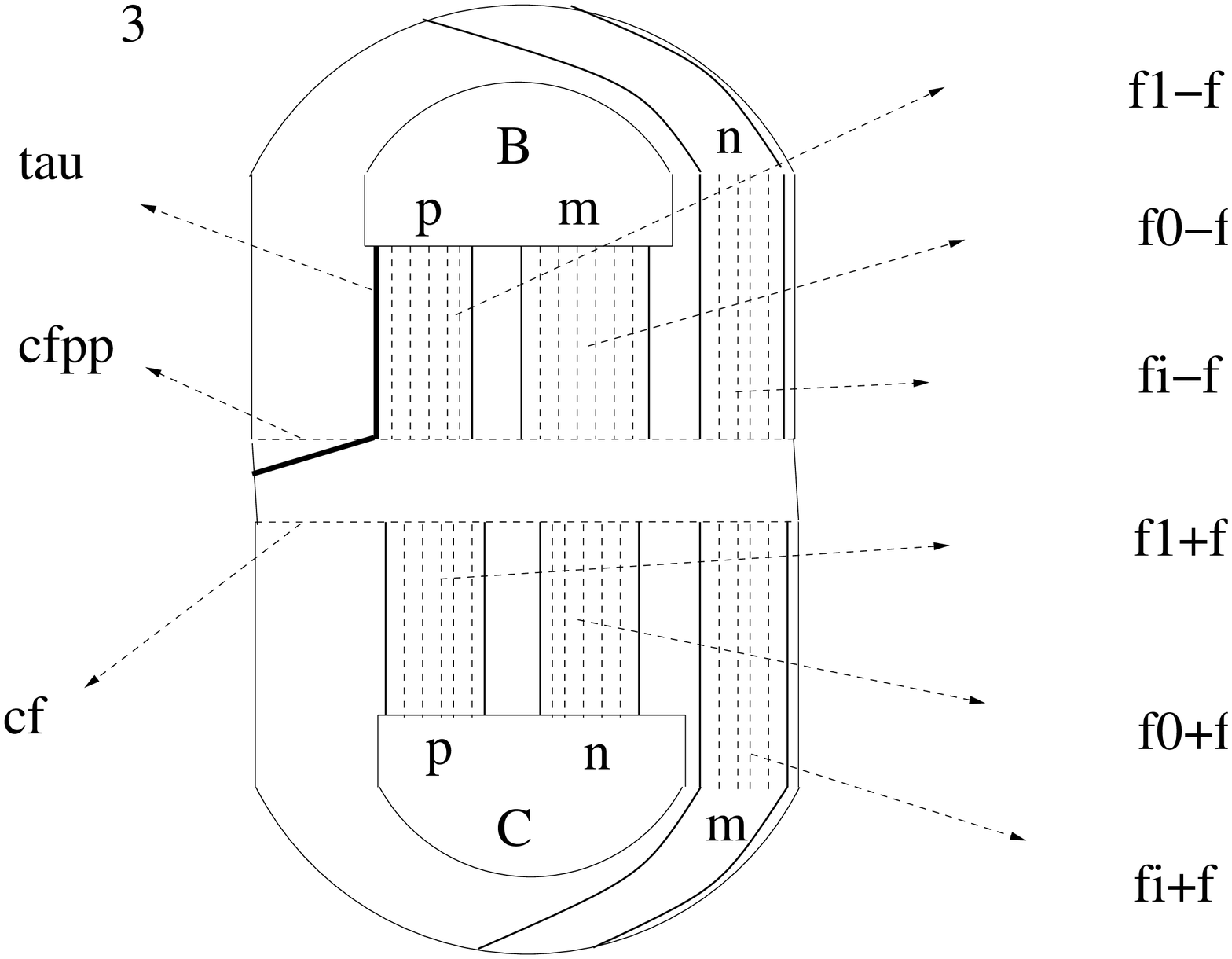}
\end{center}
\caption{}\label{fig:III}
\end{figure}

In the below figures, intersections of two reducing spheres $R$, $R^\prime$ with $T$ are shown. It is easy to see that $R\in v_S$ and  $R^\prime\in v_{S^\prime}$. 

\newpage
\begin{description}
\item[III.A]  If   $\{g_{ij}\}\subseteq\{h_{ij}\}$  (see figure \ref{fig:IIIA}): Let $x=|\{h_{ij}\}\cap \{k_{ij}\}|/2$. Then by an argument similar to the proof of the case I.A.1. we get $2(n+m+p)=v_P\cdot v_Q>v_{R^\prime}\cdot v_Q=2(m+p-n)>v_R\cdot v_Q=2(m+p-n-2x)$, \ $v_{R} \cdot v_{R^\prime}=4$ and $v_{\beta^i (R)} \cdot v_Q,\ v_{\beta^i(R^\prime)}\cdot v_Q > v_P\cdot v_Q$ for $i\neq 0$.
\begin{figure}[h]
\begin{center}
\psfrag{cr}[][][1.5]{$c_R$}
\psfrag{crp}[][][1.5]{$c_{R^\prime}$}
\psfrag{B}[][][1]{$B$}
\psfrag{C}[][][1]{$C$}
\psfrag{3}[][][1.5]{$P_f$}
\includegraphics[height=6cm, width=6cm]{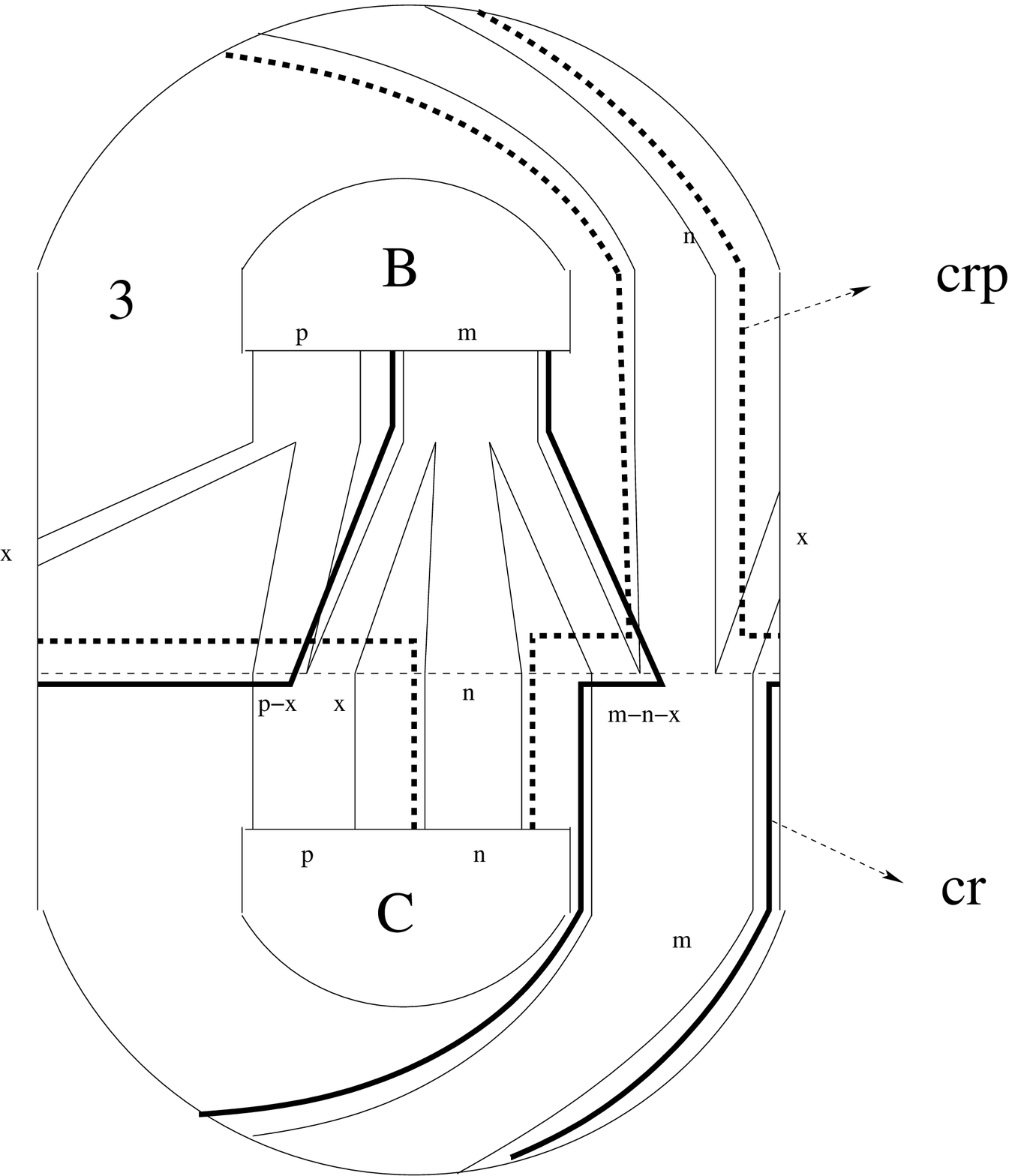}
\end{center}
\caption{}\label{fig:IIIA}
\end{figure}
\item[III.B]  If  $\{g_{ij}\}\cap\{h_{ij}\}\neq \emptyset$,  $\{g_{ij}\}\cap\{f_{ij}\}\neq \emptyset$, $\{e_{ij}\}\cap\{h_{ij}\}=\emptyset$ (see figure \ref{fig:IIIB}):  Let $x=|\{k_{ij}\}\cap \{h_{ij}\}|/2$. Then by an argument similar to the proof of the case I.A.1. we get $2(n+m+p)=v_P\cdot v_Q>v_{R^\prime}\cdot v_Q=2(p+n-m+2x)>v_R\cdot v_Q=2(p+n-m)$, \  $v_R\cdot v_{R\prime}=4$ and   $v_{\beta^i (R)} \cdot v_Q,\ v_{\beta^i(R^\prime)}\cdot v_Q > v_P\cdot v_Q$ for $i\neq 0$.
\begin{figure}[h]
\begin{center}
\psfrag{cr}[][][1.5]{$c_R$}
\psfrag{crp}[][][1.5]{$c_{R^\prime}$}
\psfrag{B}[][][1]{$B$}
\psfrag{C}[][][1]{$C$}
\psfrag{3}[][][1.5]{$P_f$}
\includegraphics[height=6cm, width=6cm]{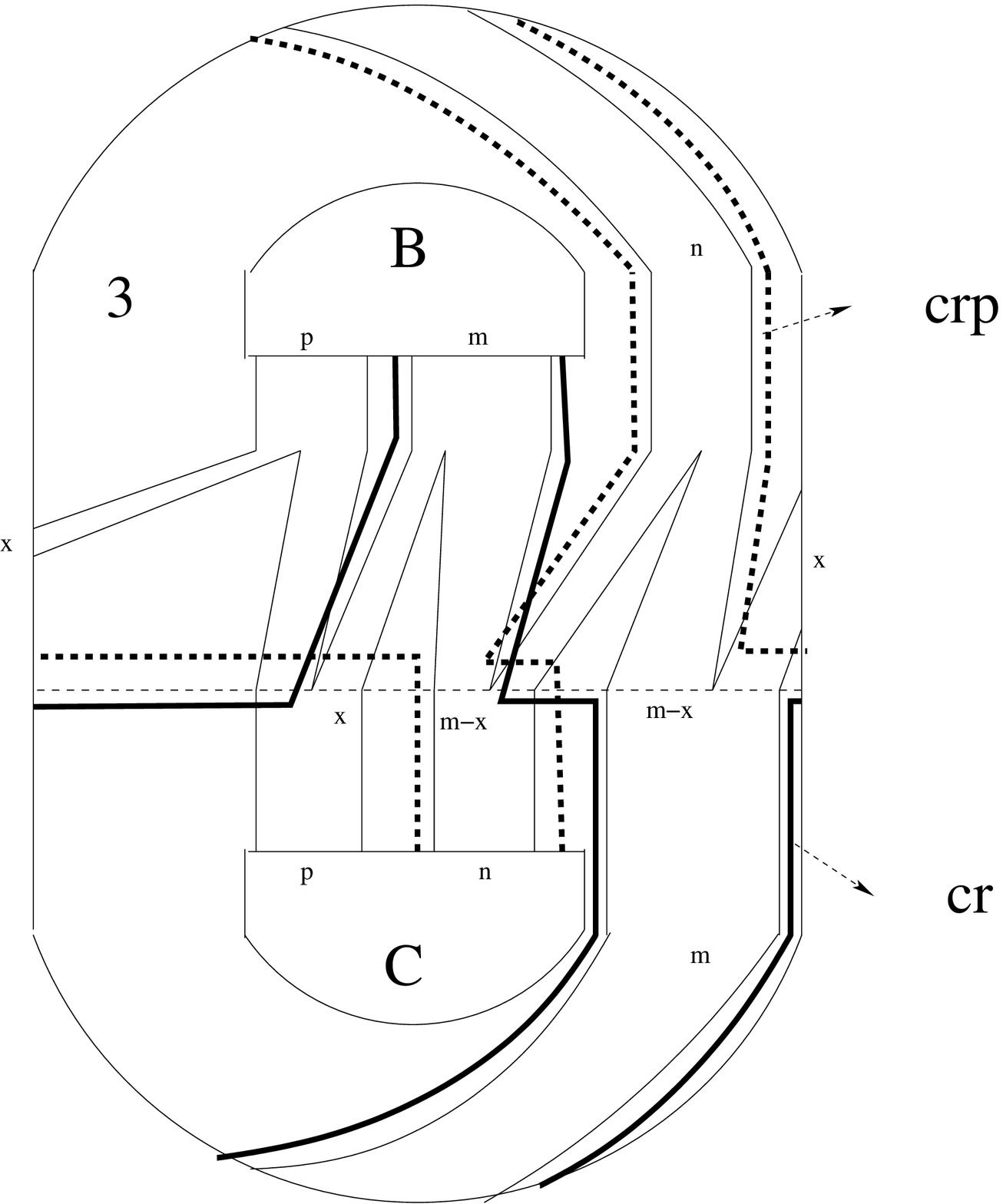}
\end{center}
\caption{}\label{fig:IIIB}
\end{figure}

\newpage
\item[III.C]  If $\{g_{ij}\}\cap\{h_{ij}\}\neq \emptyset$, $\{g_{ij}\}\cap\{f_{ij}\}\neq \emptyset$, $\{e_{ij}\}\cap\{h_{ij}\}\neq \emptyset$ (see figure \ref{fig:IIIC}):  Let $x=|\{f_{ij}\}\cap \{g_{ij}\}|/2$. Then by an argument similar to the proof of the case I.A.1. we get  $2(n+m+p)=v_P\cdot v_Q>v_{R^\prime}\cdot v_Q=2(m-n+2x+p)>v_R\cdot v_Q=2(m-n-p+2x)$, \  $v_R\cdot v_{R\prime}=4$ and  $v_{\beta^i (R)} \cdot v_Q,\ v_{\beta^i(R^\prime)}\cdot v_Q > v_P\cdot v_Q$ for $i\neq 0$.
\begin{figure}[h]
\begin{center}
\psfrag{cr}[][][1.5]{$c_R$}
\psfrag{crp}[][][1.5]{$c_{R^\prime}$}
\psfrag{E}[][][1]{$E$}
\psfrag{F}[][][1]{$F$}
\psfrag{B}[][][1]{$B$}
\psfrag{C}[][][1]{$C$}
\psfrag{3}[][][1.5]{$P_f$}
\includegraphics[height=5.7cm, width=6cm]{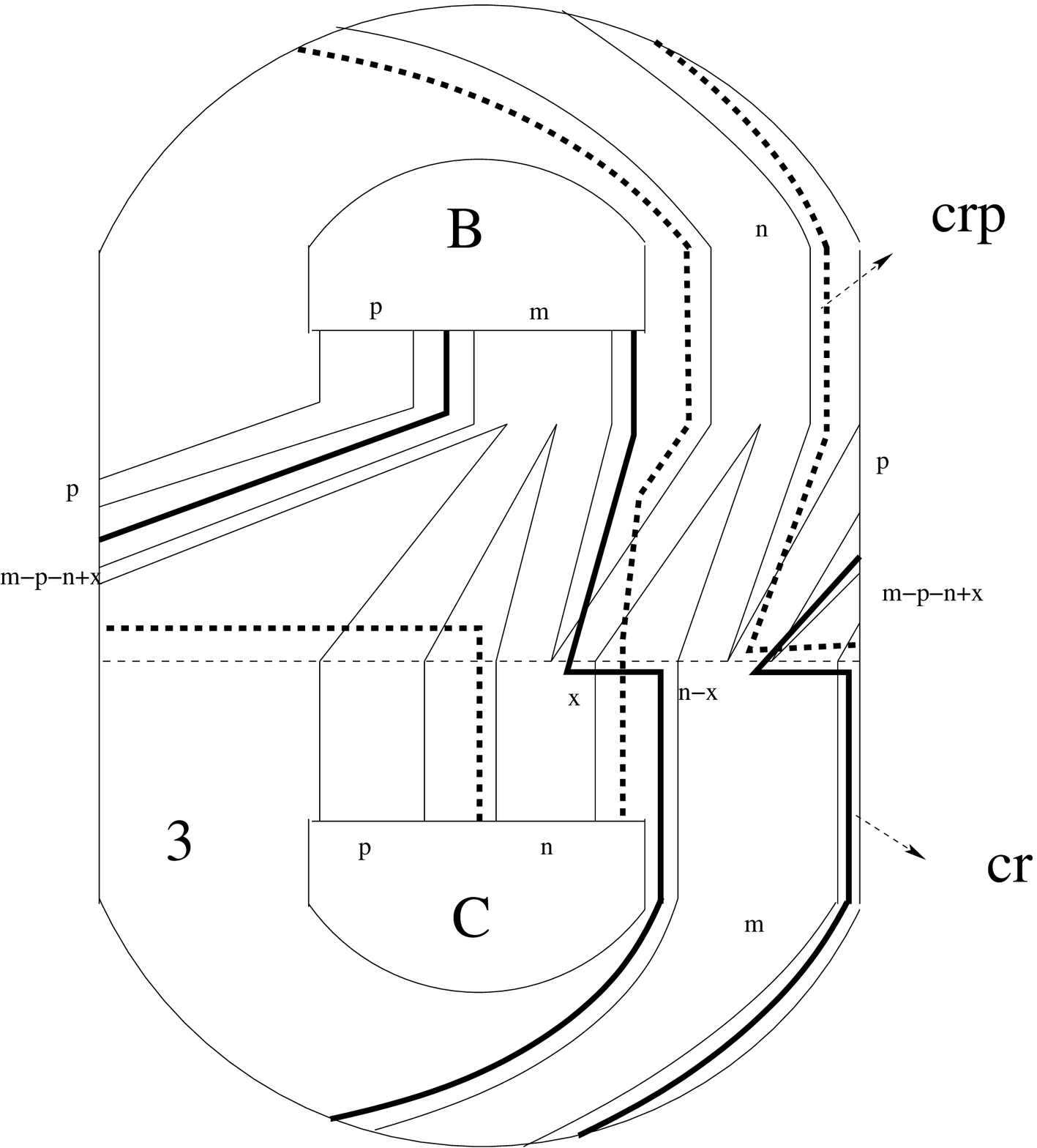}
\end{center}
\caption{}\label{fig:IIIC}
\end{figure}
\item[III.D]  If  $\{g_{ij}\}\cap\{f_{ij}\}\neq \emptyset$, $\{g_{ij}\}\cap\{l_{ij}\}\neq \emptyset$,  $\{e_{ij}\}\cap\{l_{ij}\}\neq \emptyset$,  $\{e_{ij}\}\cap\{h_{ij}\}= \emptyset$ (see figure \ref{fig:IIID}):  Let $x=|\{g_{ij}\}\cap \{l_{ij}\}|/2$. Then by an argument similar to the proof of the case I.A.1. we get $2(n+m+p)=v_P\cdot v_Q>v_{R^\prime}\cdot v_Q=2(p+n+m-2x)>v_R\cdot v_Q=2(p+n-m)$, \  $v_R\cdot v_{R\prime}=4$ and $v_{\beta^i (R)} \cdot v_Q,\ v_{\beta^i(R^\prime)}\cdot v_Q > v_P\cdot v_Q$ for $i\neq 0$.
\begin{figure}[h]
\begin{center}
\psfrag{cr}[][][1.5]{$c_R$}
\psfrag{crp}[][][1.5]{$c_{R^\prime}$}
\psfrag{B}[][][1]{$B$}
\psfrag{C}[][][1]{$C$}
\psfrag{3}[][][1.5]{$P_f$}
\includegraphics[height=5.7cm, width=6cm]{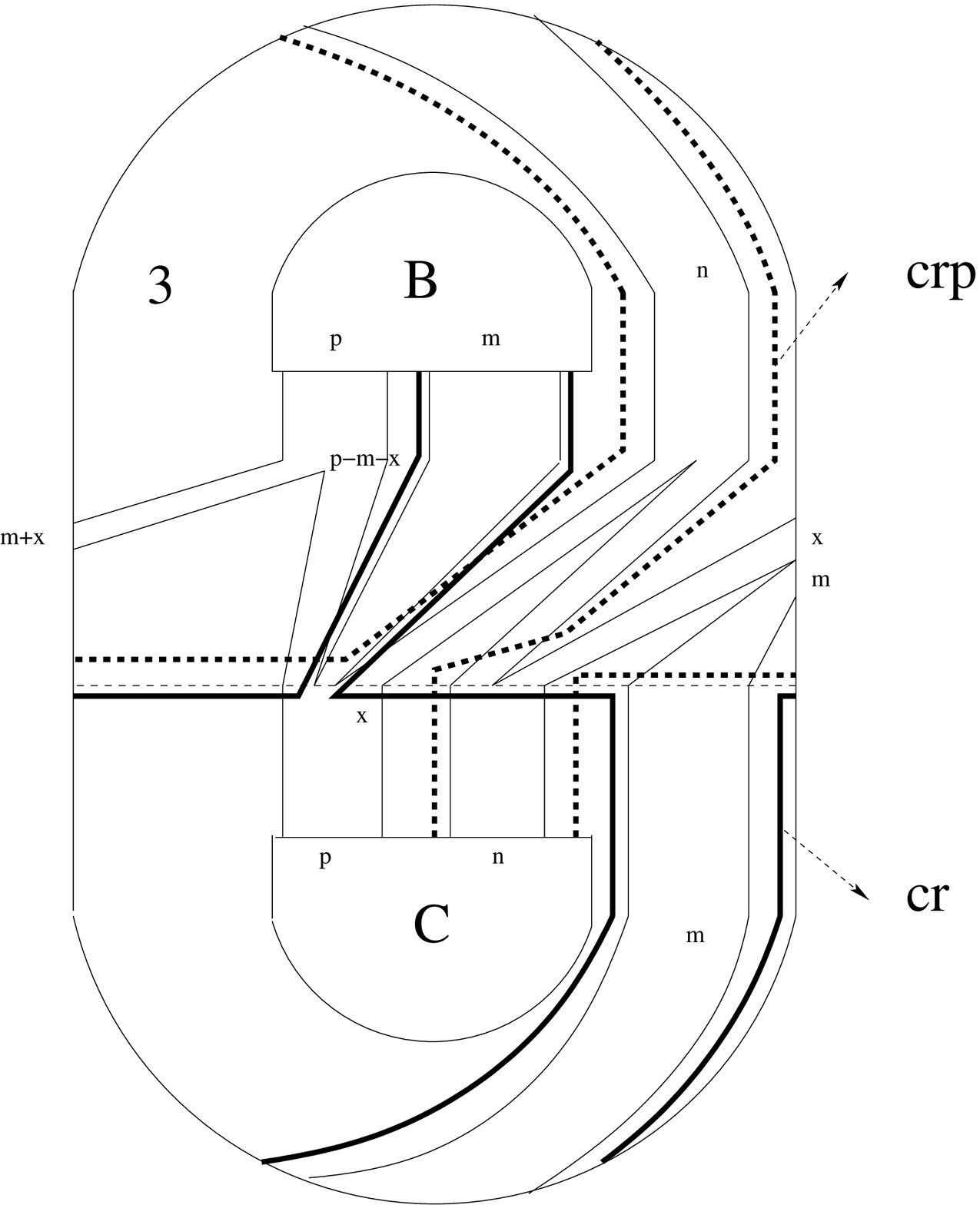}
\end{center}
\caption{}\label{fig:IIID}
\end{figure}

\item[III.E]  If $\{g_{ij}\}\cap\{f_{ij}\}\neq \emptyset$,  $\{g_{ij}\}\cap\{l_{ij}\}\neq \emptyset$, $\{e_{ij}\}\cap\{l_{ij}\}\neq \emptyset$,  $\{e_{ij}\}\cap\{h_{ij}\}\neq \emptyset$ (see figure \ref{fig:IIIE}): Let $x=|\{g_{ij}\}\cap \{l_{ij}\}|/2$. Then  by an argument similar to the proof of the case I.A.1. we get  $v_{R^\prime}\cdot v_Q=2(m+n+p-2x)$,   $v_{R}\cdot v_Q=2(m+n-p+2x)$ and $v_{\beta^i (R)} \cdot v_Q,\ v_{\beta^i(R^\prime)}\cdot v_Q > 2(n+m+p)$ for $i\neq 0$.   So $v_{R^\prime}\cdot v_Q=v_{R}\cdot v_Q$ if and only if   $p=2x$. If $p$ is equal to $2x$ then by an argument given in the proof of Lemma \ref{lemmaC} we can show that $c_Q$ does not bound a disc in $V$. Therefore either   $v_{R^\prime}\cdot v_Q>v_{R}\cdot v_Q$ or  $v_{R^\prime}\cdot v_Q<v_{R}\cdot v_Q$.  Notice that $v_R\cdot v_{R\prime}=4$. 
\begin{figure}[h]
\begin{center}
\psfrag{cr}[][][1.5]{$c_R$}
\psfrag{crp}[][][1.5]{$c_{R^\prime}$}
\psfrag{B}[][][1]{$B$}
\psfrag{C}[][][1]{$C$}
\psfrag{3}[][][1.5]{$P_f$}
\includegraphics[height=7cm, width=7cm]{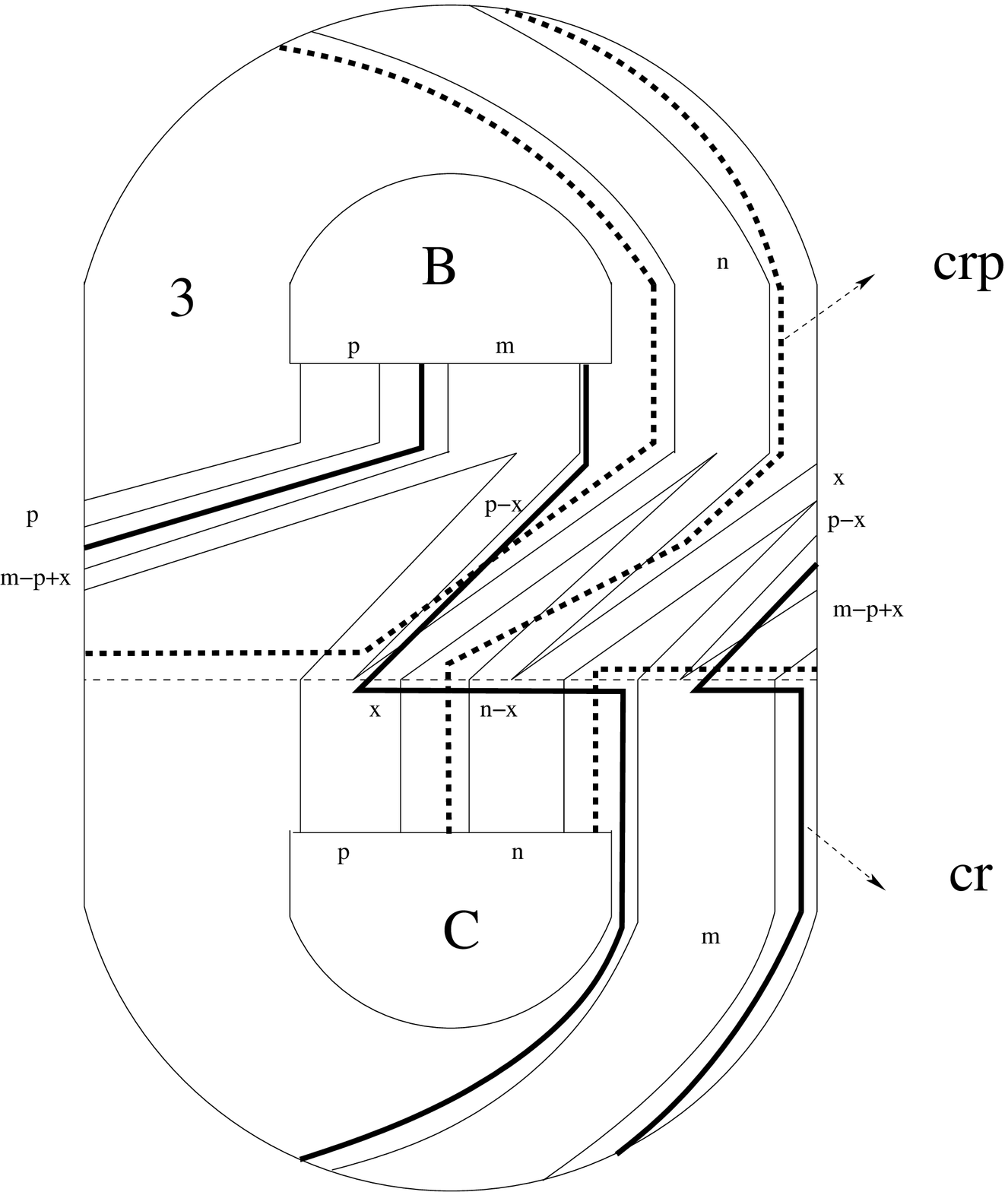}
\end{center}
\caption{}\label{fig:IIIE}
\end{figure}
\item[III.F]  If  $\{g_{ij}\}\cap\{f_{ij}\}\neq \emptyset$,  $\{g_{ij}\}\cap\{l_{ij}\}\neq \emptyset$,  $\{g_{ij}\}\cap\{h_{ij}\}\neq \emptyset$  (see figure \ref{fig:IIIF}):  Let $x=|\{g_{ij}\}\cap \{f_{ij}\}|/2$ then by an argument similar to the proof of the case I.A.1.  we get $2(n+m+p)=v_P\cdot v_Q>v_{R}\cdot v_Q=2(m+x+3p-n+x)>v_{R^\prime}\cdot v_Q=2(m+x+p-n+x)$, \  $v_R\cdot v_{R\prime}=4$ and $v_{\beta^i (R)} \cdot v_Q,\ v_{\beta^i(R^\prime)}\cdot v_Q > v_P\cdot v_Q$ for $i\neq 0$.
\newpage
\begin{figure}[h]
\begin{center}
\psfrag{cr}[][][1.5]{$c_R$}
\psfrag{crp}[][][1.5]{$c_{R^\prime}$}
\psfrag{B}[][][1]{$B$}
\psfrag{C}[][][1]{$C$}
\psfrag{3}[][][1.5]{$P_f$}
\includegraphics[height=6cm, width=6cm]{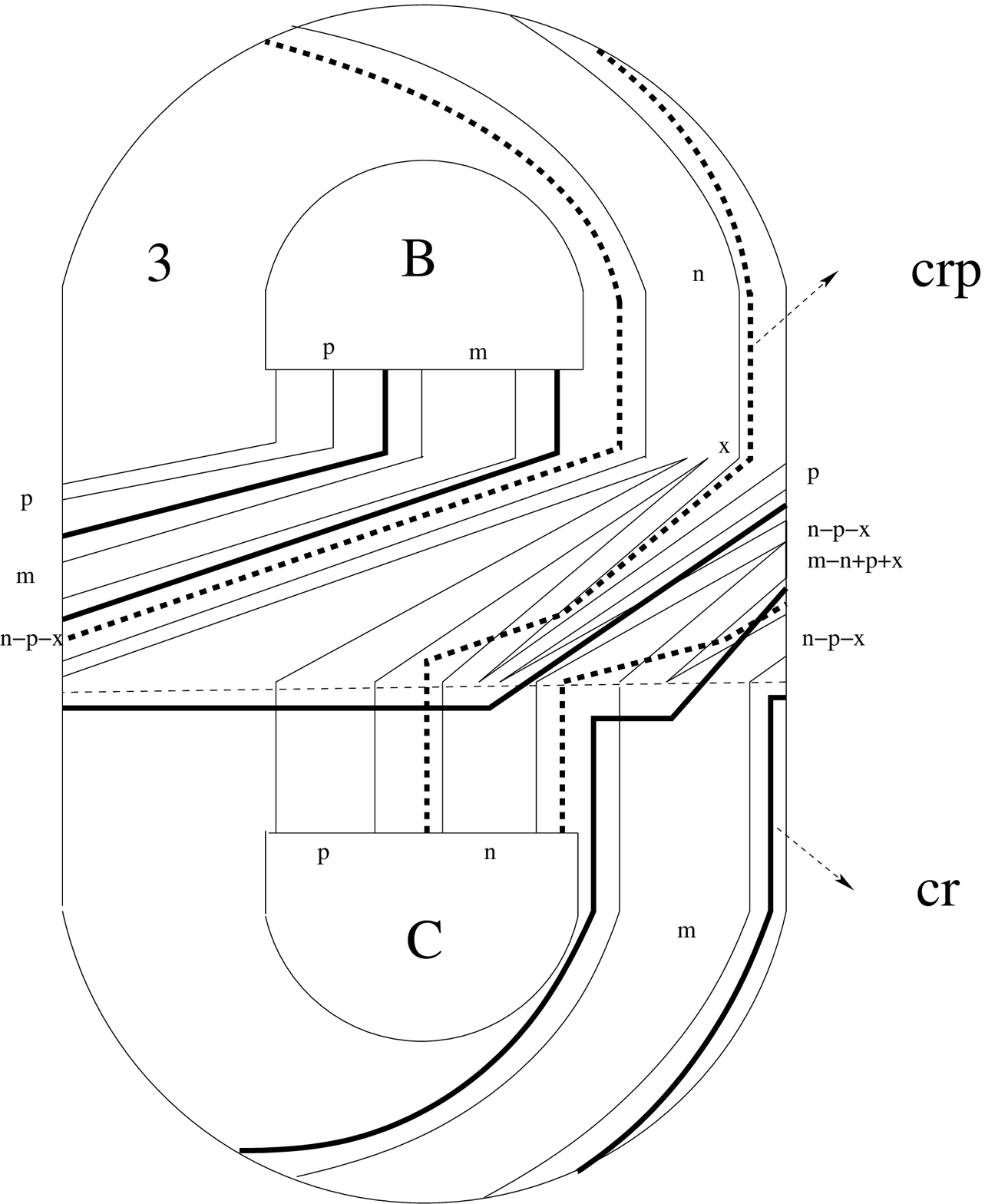}
\end{center}
\caption{}\label{fig:IIIF}
\end{figure}
\item[III.G]  If  $\{e_{ij},g_{ij}\}\subseteq\{l_{ij}\}$  (see figure \ref{fig:IIIG}):  Let $x=|\{k_{ij}\}\cap \{l_{ij}\}|/2$  then by an argument similar to the proof of the  case I.A.1. we get $v_P\cdot v_Q>v_{R^\prime}\cdot v_Q=2(p+m-n)>v_{R}\cdot v_Q=2(p-m+n)$, \   $v_R\cdot v_{R\prime}=4$ and  $v_{\beta^i (R)} \cdot v_Q,\ v_{\beta^i(R^\prime)}\cdot v_Q > v_P\cdot v_Q$ for $i\neq 0$.
\begin{figure}[h]
\begin{center}
\psfrag{cr}[][][1.5]{$c_R$}
\psfrag{crp}[][][1.5]{$c_{R^\prime}$}
\psfrag{B}[][][1]{$B$}
\psfrag{C}[][][1]{$C$}
\psfrag{3}[][][1.5]{$P_f$}
\includegraphics[height=6cm, width=6cm]{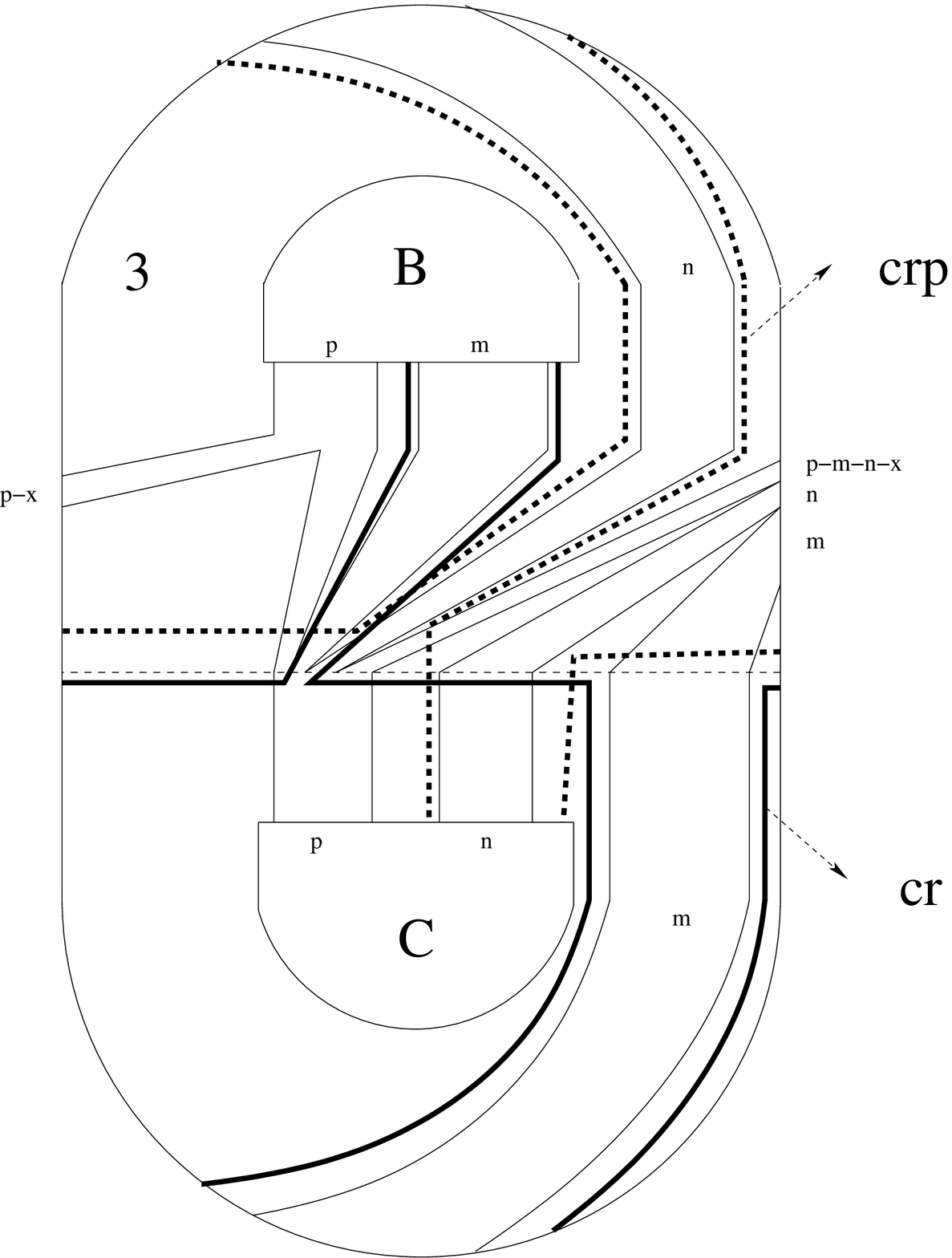}
\end{center}
\caption{}\label{fig:IIIG}
\end{figure}
\item[III.H]  If  $\{g_{ij}\}\subseteq\{l_{ij}\}$,   $\{e_{ij}\}\cap\{h_{ij}\}\neq \emptyset$: This case is eliminated  by an argument given in proof of Lemma \ref{lemmaC} (the curve $c_Q$ does not  bound a disc in $V$).
\item[III.I]  If   $\{g_{ij}\}\cap\{l_{ij}\}\neq \emptyset$,   $\{g_{ij}\}\cap\{h_{ij}\}\neq \emptyset$: After applying  $\beta^{-1}$ to $c_Q$ we can assume that $c_Q$ is as in figure \ref{fig:IIII}.  Let $x=|\{k_{ij}\}\cap \{l_{ij}\}|/2$  then by an argument similar to the proof of the case I.A.1.  we have $2(n+m+p)=v_P\cdot v_Q>v_{R}\cdot v_Q=2(m-n+3p-2x)  >v_{R^\prime}\cdot v_Q=2(m-n+p)$, \  $v_R\cdot v_{R\prime}=4$ and  $v_{\beta^i (R)} \cdot v_Q,\ v_{\beta^i(R^\prime)}\cdot v_Q > v_P\cdot v_Q$ for $i\neq 0$.
\begin{figure}[h]
\begin{center}
\psfrag{cr}[][][1.5]{$c_R$}
\psfrag{crp}[][][1.5]{$c_{R^\prime}$}
\psfrag{B}[][][1]{$B$}
\psfrag{C}[][][1]{$C$}
\psfrag{3}[][][1.5]{$P_f$}
\includegraphics[height=5.5cm, width=5.5cm]{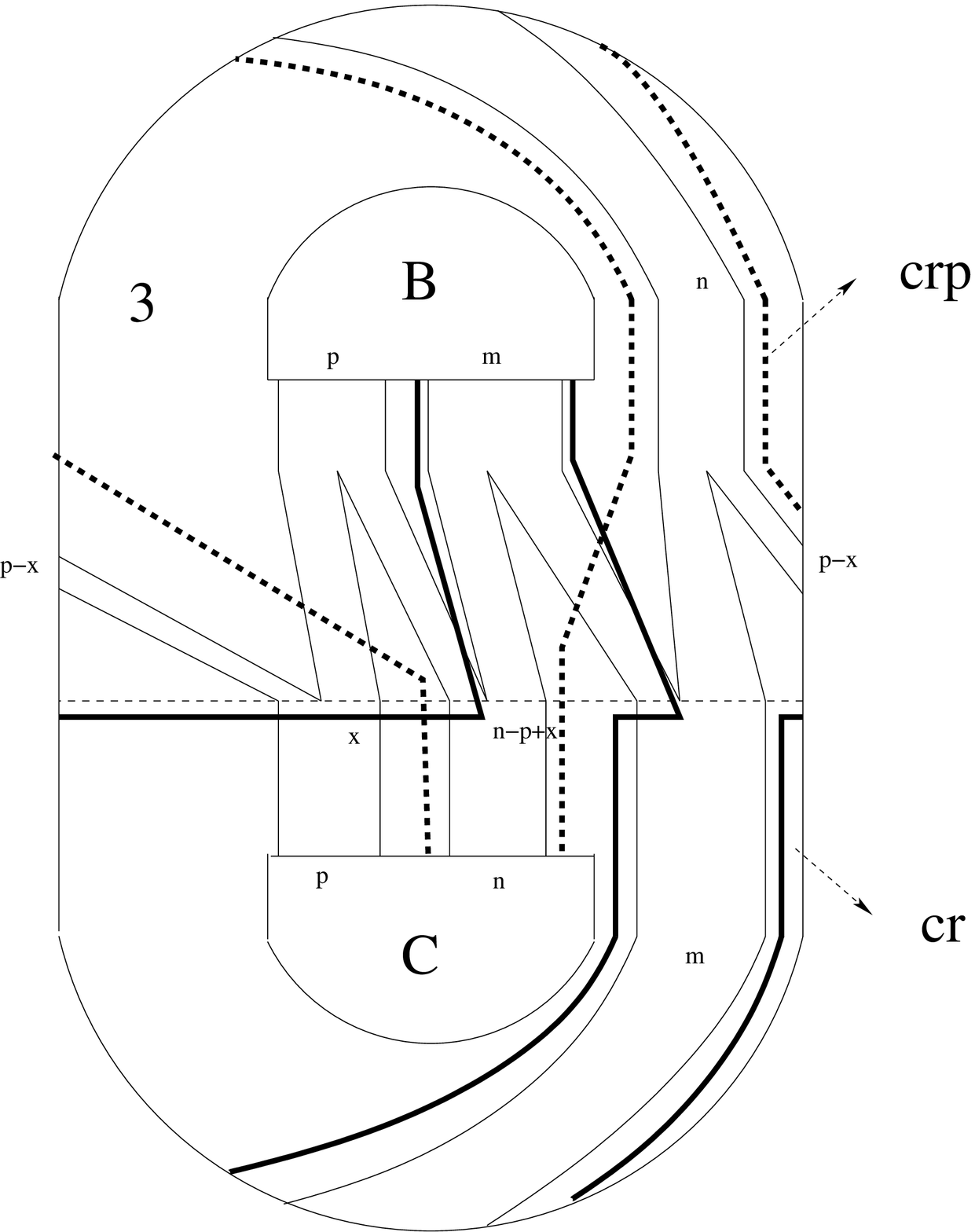}
\end{center}
\caption{}\label{fig:IIII}
\end{figure}
\end{description}
\end{proof}

\section{A presentation for $\mathcal{H}_2$}
We will first prove Theorem \ref{tree}. Then by using Bass-Serre theory we will prove Theorem \ref{presentation}.

{\it Proof of Theorem \ref{tree}.}

\begin{proof}
Suppose that  $\tilde{\Gamma}$ is not a tree. Then there is  a nontrivial loop in $\tilde{\Gamma}$. For any loop  $\xi$ in $\tilde{\Gamma}$ let $NV(\xi)$ denote the number of vertices  of $\xi$. Then $\alpha_0=\min \{ NV(\xi)\ |$ $\xi$ is a nontrivial  loop in $\tilde{\Gamma}$  $\} > 0$. Since each edge of $\Gamma$ lies  on a single $2$-simplex $\alpha_0 \geq 8$. Let $\xi_0$ be a  nontrivial loop in  $\tilde{\Gamma}$ such that $NV(\xi_0)=\alpha_0$. Since $\xi_0$ is of minimal length all its vertices are distinct. Let $v_0$ be any vertex of $\xi_0$, and let $v_0$,  $v_1$, $v_2$, $v_3$,..., $v_{\alpha_0-1}$  be the consecutive  vertices of $\xi_0$.  We may suppose that $v_0\in\Gamma$. Then $v_0$, $v_2$, $v_4$,..., $v_{\alpha_0-2}$ are  vertices of $\Gamma$, and $v_k\cdot v_{k+2}=v_{\alpha_0-2}\cdot v_0=4$ for $k \in \{0,2,4,...,\alpha_0-4\}$. 

{\bf Claim.}  $v_k\cdot v_0 < v_{k+2}\cdot v_0$ for $k \in \{0,2,4,...,\alpha_0-4\}$.

{\bf Proof of claim.} The proof will be by induction on the index $k$. If $k=0$ then $v_0\cdot v_0=0 < v_{2}\cdot v_0=4$. Assume  $v_k\cdot v_0 < v_{k+2}\cdot v_0$ for $k \in \{0,2,...,\alpha_0-6\}$. If $v_{k+4}\cdot v_0\leq v_{k+2}\cdot v_0$ then $v_{k}\cdot v_{k+4}=4$ by Proposition \ref{backbone}. Since $v_{k}\cdot v_{k+2}=v_{k+2}\cdot v_{k+4}=4$, the vertices  $v_{k}$, $v_{k+2}$, $v_{k+4}$ form a 2-simplex $\triangle$ in  $\Gamma$. Then  we get a loop $\xi$ in $\tilde{\Gamma}$ with vertices  $v_0$,  $v_1$,..., $v_{k}$, $u$, $v_{k+4}$, $v_{k+5}$, ... , $v_{\alpha_0-2}$, $v_{\alpha_0-1}$  where $u$ is the barycenter of $\triangle$.  This contradicts the minimality of  $\alpha_0$.

By the above claim,  we get $v_0\cdot v_{\alpha_0-4}<v_0\cdot v_{\alpha_0-2}$. But $4<v_0\cdot v_{\alpha_0-4}$ and $v_0\cdot v_{\alpha_0-2}=4$, a contradiction.
\end{proof}

{\it Proof of Theorem \ref{presentation}.}

\begin{proof}
Let $v_M$ be a vertex of  $\tilde{\Gamma}$ corresponding to the barycenter of the  $2$-simplex whose vertices are $v_P$, $v_{\delta(P)}$ and $v_{\delta^2(P)}$. Let $E$ be the edge of $\tilde{\Gamma}$ whose vertices are $v_P$ and $v_M$. Let  $H_P$ be the subgroup of $\mathcal{H}_2$ generated by the elements that stabilize $v_P$. Let $H_M$ be the subgroup of $\mathcal{H}_2$ generated by the elements that preserve $v_M$.  Let $H_E$ be the group of elements of $\mathcal{H}_2$ that stabilize the edge $E$. 

\begin{itemize}

\item Scharlemann in \cite[Lemma 2]{Sc} gives the following  presentation for $H_P$:

$H_P=\ <[\alpha],[\beta],[\gamma]\ |\ [\alpha]^2=[\gamma]^2= [\alpha\gamma]^2=[\alpha\beta\alpha\beta^{-1}]=1,[\gamma\beta\gamma]=[\alpha\beta]>\ \cong (\mathbb{Z} \oplus \mathbb{Z}_2) \rtimes \mathbb{Z}_2$.

\item The subgroup  $H_M$ fixes the set  $\{v_P,v_{\delta(P)},v_{\delta^2(P)}\}$. Therefore

$H_M=<[\delta],[\alpha],[\gamma]\ |\ [\delta]^3=[\alpha]^2=[\gamma]^2= [\alpha\delta\alpha^{-1}\delta^{-1}]=[\alpha\gamma]^2=1,\ \ [\delta] =[\gamma\delta^2\gamma]> \ \
\cong  (\mathbb{Z}_3 \rtimes \mathbb{Z}_2) \oplus \mathbb{Z}_2$.

\item An element $h$ of $\mathcal{H}_2$ fixes the sets $\{v_P\}$ and $\{v_{\delta(P)},v_{\delta^2(P)}\}$ if and only if $h\in H_E$. Hence

$H_E=<[\alpha],[\gamma]|\ \ [\alpha]^2 =[\gamma]^2 = [\alpha\gamma]^2=1>\cong \mathbb{Z}_2 \oplus \mathbb{Z}_2$.

\end{itemize}

The action of $\mathcal{H}_2$ on the 2-complex $\Gamma$ induces an action of $\mathcal{H}_2$ on the tree $\tilde{\Gamma}$. The subgroups $H_P$,  $H_M$ are the isotropy subgroups of $\mathcal{H}_2$ fixing the vertices $v_P,v_M$ respectively. By the standard Bass-Serre theory \cite{S} the group $\mathcal{H}_2$ is thus a free product of the subgroups $H_P$ and $H_M$ amalgamated over the
subgroup $H_E$.

\[
\mathcal{H}_2 \cong  H_P \underset{H_E}{\ast} H_M \cong  (\mathbb{Z} \oplus \mathbb{Z}_2) \rtimes \mathbb{Z}_2 \underset{\mathbb{Z}_2 \oplus \mathbb{Z}_2}{\ast} (\mathbb{Z}_3 \rtimes \mathbb{Z}_2) \oplus \mathbb{Z}_2
\]

\end{proof}

\end{document}